\documentclass[12pt,a4paper]{amsart}

\usepackage{latexsym}

\usepackage{amssymb}
\usepackage{amsthm}
\usepackage{amsfonts}
\usepackage{amsmath}
\usepackage{amstext}
\usepackage{amscd}
\usepackage[dvips]{graphics}
\usepackage{epsfig}

\numberwithin{equation}{section}

\theoremstyle{plain}%default
\newtheorem{thm}{Theorem}[section] 
\newtheorem{prop}[thm]{Proposition}
\newtheorem{cor}[thm]{Corollary}
\newtheorem{lem}[thm]{Lemma}

\newtheorem{theorem*}{Theorem}[]

\theoremstyle{definition}
\newtheorem{defn}[thm]{Definition}

\newtheorem{example}[thm]{Example}

\theoremstyle{remark}
\newtheorem{rem}[thm]{Remark}

\newcommand{\C}{\mathbb{C}}
\newcommand{\N}{\mathbb{N}}
\newcommand{\R}{\mathbb{R}}

\newcommand{\Q}{\mathbb{Q}}
\newcommand{\Z}{\mathbb{Z}}

\newcommand{\projc}{\mathbb{CP}}
\newcommand{\projr}{\mathbb{RP}}
\newcommand{\compl}[1]{{#1}_\mathbb{C}}

\newcommand{\inv}{^{-1}}

\newcommand{\wR}{\widetilde {\mathbb{R}^2}}

\DeclareMathOperator{\mult}{mult\,}
\DeclareMathOperator{\ord}{ord\,}
\DeclareMathOperator{\sgn}{sign\,}

\DeclareMathOperator{\dist}{dist\,}

\title[Blow-analytic equivalence] {Blow-analytic equivalence of two variable 
real analytic function germs}

\author{Satoshi Koike \& Adam Parusi\'nski}

\address {Department of Mathematics, Hyogo University
of Teacher Education, 942-1 Shimokume, Kato, 
Hyogo 673-1494, Japan}

\email {koike@hyogo-u.ac.jp}

\address {Laboratoire Angevin de Recherche en Math\'ematiques, UMR 6093 
du CNRS,  Universit\'e d'Angers,
   2, bd Lavoisier, 49045 Angers cedex, France}

\email{adam.parusinski@univ-angers.fr}

 
\thanks{This research was partially supported by the Grant-in-Aid 
for Scientific Research (No. 18540084) of Ministry of Education, 
Science and Culture of Japan. }

\begin{document}


%%%%%%%%%%%%%%%%%%%%%%%%%%%%%%%%%%%%%%%%%%%%%%%%%%%%%%%%%%%%%%%%%%%%

%\vspace{0.2 truecm}  

\begin{abstract} 
Blow-analytic equivalence is a notion for real analytic function
germs, 
introduced by Tzee-Char Kuo in order to develop real analytic 
equisingularity theory. In this paper we give complete 
 characterisations of blow-analytic equivalence in the two 
dimensional case:  in terms of  the real tree 
model for the arrangement of real parts of  Newton-Puiseux roots and their Puiseux pairs, and 
in terms  of  minimal resolutions.  
 These characterisations show that in the
two dimensional case the blow-analytic equivalence is a natural 
analogue of topological equivalence of complex analytic function
germs. 
Moreover, we show that in the two-dimensional case the blow-analytic equivalence can be 
made cascade, and hence satisfies several  geometric properties.   It preserves, for instance,  
the contact orders of real analytic arcs.  

In the 
general $n$-dimensional case, we show that a singular real modification  satisfies 
the arc-lifting property.  
\end{abstract}

\maketitle

%%%%%%%%%%%%%%%%%%%%%%%%%%%%%%%%%%%%%%%%%%%%%%%%%%%%%%%%%%%%%%%%%%%%

%\section{Introduction}
%\label{intro}
%\medskip

\vspace{0.2 truecm}
A classical result of   Burau \cite{burau} and Zariski \cite{zariski} shows the embedded  topological type of 
a plane curve singularity  $(X , 0)\subset  (\C^2,0)$ is determined by the Puiseux pairs of each 
irreducible component   and the intersection numbers of any pairs of 
distinct components.  It can be shown, cf. \cite {parusinski}, that  the topological 
type of function germs $f: (\C^2,0) \to (\C,0)$ is
completely characterised, also in the non-reduced case $f=\prod f_i^{d_i}$,  by the embedded topological type of its zero set and the multiplicities $d_i$ of its irreducible components.  

In this paper we give a real 
analytic counterpart of these results and show that the  two variable 
version of blow-analytic equivalence of Kuo is classified by  invariants similar to Puiseux pairs, multiplicities of irreducible components,  and intersection numbers. 
 Moreover we show several natural geometric properties  of this equivalence, answering previously posed questions.  
In the main result of this paper we a give complete characterisation of blow-analytic
equivalence classes of two variable real analytic function germs. 

\begin{thm}\label{allequivalent}
Let $f:(\R^2,0)\to (\R,0)$ and $g:(\R^2,0)\to (\R,0)$ be real analytic 
function germs.  Then the following conditions are equivalent:
\begin{enumerate}
\item %[] 
%\item[{ }(a)]
$f$ and $g$ are blow-analytically equivalent.
%\item[{ }(b)]
%\item $f$ and $g$ are cascade-blow-analytically equivalent.
\item %[{ }(c)]
$f$ and $g$ have weakly isomorphic minimal resolution spaces.
\item %[{ }(d)]
The real tree models of $f$ and $g$ are isomorphic. 
%\item[{ } (e)] 
%The family of real singular modifications is stable by taking the 
%strict transform by local blowings up:  
%$f$ and $g$ can be connected by an equi-resolvable deformation.  
\end{enumerate}
Moreover if $f$ and $g$ are blow-analytically equivalent then 
they are equivalent by a cascade-blow-analytic homeomorphism. 
\end{thm}

Theorem \ref{allequivalent} can be stated in both the oriented 
and non-oriented case, see section \ref{EndofProof} below.  
By a weak isomorphism of resolution spaces we mean a 
 homeomorphism that preseves the basic numerical data of resolutions, 
see subsection \ref{constructing}.  
The real tree model is a counterpart of Kuo and Lu's tree model \cite{kuolu}, 
a combinatorial object that encodes the numerical data 
given by the contact orders between the Newton-Puiseux roots of $f$ in the complex case.  
Cascade-blow-analytic homeomorphisms satisfy many geometric and analytic properties: 
they lift to the resolution spaces of $f$ and $g$,  they  preserve the intersection numbers between real analytic arcs and their Puiseux exponents.  
%We will give precise definitions and statements later.  

In the general  ($n$-dimensional) case 
the notion of blow-analytic equivalence is very technical but we need to recall its definition 
and  main properties.

\subsection{Blow-analytic equivalence. } \label{b.a.e.}
In a search for a "right" equivalence relation of real analytic function
germs, that could play a similar role to the topological equivalence 
in the complex analytic set-up, at the end of 1970 Tzee-Char Kuo proposed 
the notion of blow-analytic equivalence 
\cite{kuo2, kuo3, kuo4, kuoward, kuo5, kuo6, kuo7, kuo8}).
In \cite{kuo8}, Kuo proved that blow-analytic equivalence
is an equivalence relation and established the  local finiteness
of blow-analytic types for analytic families of
real analytic function-germs with isolated singularities.

 We say that a homeomorphism germ $h :
(\R^n,0) \to (\R^n,0)$
is  a {\em blow-analytic homeomorphism}  if there exist real modifications 
$\mu : (M,\mu^{-1}(0)) \to (\R^n,0)$, 
$\tilde \mu : (\tilde M, \tilde \mu\inv (0) )$ $ \to (\R^n,0)$
and an analytic isomorphism $\Phi : (M,\mu^{-1}(0)) \to
(\tilde M,\tilde\mu^{-1}(0))$
so that the following diagram is commutative: 
\begin{equation}\label{baequivalence}
\minCDarrowwidth 1pt
\begin{CD}
(M, {\mu} \inv (0)) @>\mu>>  (\R^n,0)  \\ 
@V\Phi VV  @ VhVV  \\
(\tilde M, \tilde {\mu} \inv (0))   
@> \tilde\mu>> (\R^n,0)  
\end{CD}
\end{equation}
We say that two real analytic function germs 
$f : (\R^n,0) \to (\R,0)$ and $g : (\R^n,0) \to (\R,0)$ 
are {\em blow-analytically equivalent} 
if there exists a blow-analytic homeomorphism  $h : (\R^n,0) \to (\R^n,0)$
such that $f = g \circ h$. 

Kuo's definition of real modification is very technical and often difficult to 
work with.  We discuss it in section \ref{realmodifications}.  
For $n=2$, we show  the following simple characterisation.

\begin{thm}\label{2variablemod}
Let $X,S$ be connected nonsingular real analytic surfaces and 
let $\sigma : X\to S$ be a proper surjective real analytic map. 
Then $\sigma$ is a real modification in the sense of Kuo if and only if it is 
a composition of point blowings-up.
\end{thm}

Several  results showing that two 
function germs are blow-analytical\-ly 
equivalent were obtained using, mostly toric, equiresolutions  by Kuo, 
Fukui-Yoshina\-ga \cite{fukuiyoshinaga}, 
Fukui-Paunescu \cite{fukuipaunescu1}, 
Abderrahmane \cite {abderrahmane1}, and others.

Invariants allowing  to distinguish different 
blow-analytic types were constructed, using the geometry of 
arc-spaces and motivic integration, 
by Fukui \cite{fukui}, Fichou \cite{fichou}, 
and in \cite{koikeparusinski}.  These constructions are based on the observation 
that blow-analytic homeomorphisms send real analytic arcs to real analytic 
arcs.  
%Since we could not find any proof of this fact in the literature 
%we show it in this paper, see Corollary \ref{preservingarcs} below.  
It follows from 
Theorem \ref{rmodproperties}  that shows that the real modifications 
satisfy the arc lifting property, compare \cite {fukui}
section 3.   

For more on the blow-analytic equivalence see the surveys  
\cite {fukuikoikekuo}, \cite {fukuipaunescu2}.

%\smallskip
\subsection{Cascade blow-analytic homeomorphisms  and their geometric properties.} 
In general, a blow-analytic homeomorphism $h$ is  not necessarily Lipschitz.  In two dimensional case, if $h$ gives 
blow analytic equivalence between 
analytic function germs, then $h$ is cascade and  satisfies many  geometric properties.  

We say that $h:(\R^2,0)\to (\R^2,0)$ is \emph{a cascade blow-analytic homeomorphism} 
if there exists a  commutative diagram 
\begin{equation}\label{cascade}
\minCDarrowwidth 1pt
\begin{CD}
(M_k, E_k) @>b_{k}>> (M_{k-1}, E_{k-1}) @>b_{k-1}>> \cdots @>b_{2}>> 
(M_1, E_1)  @>b_{1} >> (\R^2,0) \\
@V\Phi VV  @ VVh_{k-1}V @. @ VVh_{1}V @ VVhV \\
(\tilde M_k, \tilde E_k) @>\tilde b_{k}>> (\tilde M_{k-1}, \tilde E_{k-1})
@>\tilde b_{k-1}>> 
\cdots @>\tilde b_{2}>> (\tilde M_1, \tilde E_1)  
 @>\tilde b_{1} >> 
(\R^2,0) ,
\end{CD}
\end{equation}
where $b_i, \tilde b_i$ are point blowings-up, $\tilde E_i, E_i$ are
the inverse images of the origin,  $h_i$ are homeomorphisms, and
$\Phi$ is  an analytic isomorphism.  
We say that the real analytic function germs $f(x,y), g(x,y)$ 
are \emph{cascade blow-analytically equivalent} 
if there exists a cascade blow-analytic homeomorphism  $h$ such that
$f=g\circ h$.

Suppose that $f$ and $g$ are blow-analytically equivalent.  The key step in showing 
(1) $\Rightarrow$(2) of theorem 
\ref{allequivalent} is to establish the existence 
of $h_1$ in \eqref{cascade}, or equivalently, that  in \eqref{baequivalence},  
$\Phi(E_1) = \tilde E_1$, where here $E_1\subset M$, resp. 
$\tilde E_1\subset \tilde M$, denotes the strict transform of the exceptional 
divisor of first point blowing-up in $\mu$, $\tilde \mu$ resp..  
This is shown in section \ref{lift} using the combinatorial properties of dual 
graphs of real resolutions.  This also shows that blow-analytic equivalence implies the cascade one
 in two variable case.

As we show in section \ref{Cascade} the cascade blow-analytic homeomorphisms 
satisfy many important geometric properties.  They preserve  
the Puiseux characteristic sequence of real analytic arcs and, in the oriented case, 
the signs of coefficients at the Puiseux characteristic exponents.  
They preserve also the order of contact between such arcs.   
These properties are crucial for the proof of (1)$\Longleftrightarrow $(3) 
of Theorem \ref{allequivalent}.  The theory of real analytic arcs, and more precisely their demi-branches, 
is developed in section \ref{arcs}.

Kobayashi and Kuo constructed in   \cite{kobayashikuo} examples of 
blow-analytic homeomorphisms  $h: (\R^2,0) \to  (\R^2,0)$  that are
not cascade.  Their examples do not satisfy 
 $\Phi(E_1) = \tilde E_1$ and  do not preserve 
the tangency of curves and may send a smooth arc to  a singular one and vice versa.  
Such blow-analytic homeomorphism cannot give blow-analytic equivalence between 
two real analytic function germs.

\subsection{Real tree model.}
In \cite {kuolu} Kuo and Lu introduced a tree model $T(f)$ of 
a complex analytic 
function germ $f(x,y)$.  This model allows one to visualise the numerical data 
given by the contact orders between the Newton-Puiseux roots of $f$, in 
particular their Puiseux characteristic exponents. 

For $f(x,y)$ real analytic ''the real part of $T(f)$'' was proposed 
by Kurdyka and Paunescu in \cite {kurdykapaunescu}.  
In section \ref{realtreemodel} we propose a similar, but more precise, 
construction of a real tree model that determines the resolution 
process of $f$. 
Our real tree model 
is a combinatorial object that encodes the contact orders 
between the real parts of complex Newton-Puiseux roots of $f$.   
It contains the information about the signs of coefficients at 
the Puiseux characteristic exponents. 
 Thanks to the geometric properties of cascade blow-analytic homeomorphisms, 
see Theorem \ref{cascadeproperties}, the proof 
of (1) $\Longleftrightarrow $ (3) of Theorem \ref{allequivalent} is based on 
a fairly straightforward computation of the tree model of 
the blown-up singularity in terms 
of the original tree model.

%\smallskip
\subsection{Examples.} 
Abderrahmane \cite{abderrahmane2} showed that blow-analytically 
equivalent weigh\-ted homogeneous singular 
$f : (\R^2,0) \to (\R,0)$ and $g : (\R^2,0) \to (\R,0)$ 
have the same weights.  This can be also 
easily verified by Theorem \ref{allequivalent}, see Example
\ref{weighted}.

The blow-analytic (non-oriented) classification of Brieskorn 
two variable singularities $\pm x^p \pm y^q$ was obtained in 
\cite{koikeparusinski} using the Fukui invariants and the zeta
functions.    
This classification coincides with $(x,y) \to (\pm x, \pm y)$ 
classification except for the case $p$ odd 
and $q=pm$ with $m$ even. In the latter case, additionally,   
$f(x,y)= x^p-y^q$ and $g(x,y)= x^p +y^q$ are blow analytically
equivalent, cf. \cite{koikeparusinski}, but they are not analytically 
equivalent.  Moreover, $f$ and $g$ are not bi-lipschitz equivalent.  
It is shown for $f(x,y)= x^3-y^6$ and $g(x,y)= x^3 +y^6$
in  \cite{henryparusinski2}, the proof works in general.

Theorem \ref{allequivalent} allows us to complete the classification 
of Brieskorn two variable 
singularities in the oriented case.  It coincides 
with the non-oriented case except for $f(x,y)= x^p-y^q$ and 
$g(x,y)= x^p +y^q$, both $p$ and $q$ odd.  For these functions 
$f$ and $g$ are blow-analytically equivalent but not by an orientation 
preserving homeomorphism, cf. exemple \ref{Brieskornoriented}. 

The functions 
$f(x,y)= x(x^3-y^5)(x^3+y^5)$ and  $g(x,y)= x(x^3-y^5)(x^3-2 y^5)$
have 
the same Fukui invariants and zeta functions.  As follows from Theorem
\ref{allequivalent} 
they are not  blow-analytically equivalent,  see Examples
\ref{examplesimple} and \ref{tree}.

%\smallskip
\subsection{Further developement. Open questions.} 
In \cite{forthcoming} we compare various 
equivalence relations between two variable function germs.   
Namely we show that two real analytic function germs 
$f : (\R^2,0) \to (\R,0)$ and $g : (\R^2,0) \to (\R,0)$ that are 
$C^1$ equivalent have the same real tree models and consequently, by Theorem \ref{allequivalent}, 
 are blow-analytically equivalent. If we assume that $f$ and $g$  only  
bi-lipschitz equivalent, then though the contact orders between 
Newton-Puiseux roots and the real parts of complex roots of $f$ and 
$g$ are preserved, the Puiseux pairs of these roots 
can be different, and therefore their tree models and blow-analytic 
types are different.  The simplest example is given by   
$x(x^3+ y^5)$ and $x(x^3-y^5)$, these two functions are bi-lipschitz  
but not blow-analytically equivalent by an orientation preserving 
homeomorphism.      

Most of the questions answered in this paper for  functions 
of two real variables remain open in higher dimensions, in 
particular, the very questions what should be the right precise definitions 
of the blow-analytic equivalence and of the real modification.  
Our corollary \ref{preservingarcs} shows that Kuo's blow-analytic 
homeomorphisms preserve real analytic arcs in the $n$ dimensional case.  Under what assumption 
does the blow-analytic equivalence   preserve the contact 
order between real analytic arcs?  Does it satisfy metric properties
stated in 
Corollary \ref{cascadedistance} and Propostion \ref{cascadedeterminant}?  
Note that these properties would allow one to construct more
blow-analytic  invariants 
using the methods of motivic integration, as in \cite{fichou}.   

In Lemma \ref{normalcrossing} we show that a function blow-analytically equivalent to 
normal crossing is itself normalcrossing with the same exponents.  Again, we show this result 
in two variable case.  Is it true in the general case?

%{(Acknowledgements)}
%{(Notation)} We shall drop the germ notation.  

\newpage

%%%%%%%%%%%%%%%%%%%%%%%%%%%%%%%%%%%%%%%%%%%%%%%%%%%%%%%%%%%%%%%%%%%%%%%%%%%

%%%%%%%%%%%%%%%%%%%%%%%%%%%%%%%%%%%%%%%%%%%%%%%%%%%%%%%%%%%%%%%%%%%%

\bigskip
\section{Preliminaries}
\label{preliminaries}

\medskip
%%%%%%%%%%%%%%%%%%%%%%%%%%%%%%%%%%%%%%%%%%%%%%%%%%%%%%%%%%%%%%%

\subsection{Dual resolution graph.}\label{graph}

Let $f : (\R^2,0) \to (\R,0)$ be an analytic function germ.  
We call a composition of point blowings-up $\mu : M\to \R^2$ 
\emph {a resolution of $f$} if $f\circ \mu$ is normal crossings.  
Since all blowings-up have point centres there exists a unique minimal resolution of 
$f$ obtained by blowing-up only the points where the total transform of $f$ is 
not normal crossings.  

Consider the weighted oriented \emph{dual graph $S=S_{\mu,f}$ associated 
to a resolution $\mu$ of $f$}.  
Each component $E$ of the exceptional divisor $\mu\inv (0)$ corresponds to a vertex of $S$.  
For simplicity we denote this vertex also by $E$. If two such components intersects 
at a point then the corresponding vertices are joined by an edge. Each component $C$ of 
the strict transform of $f\inv (0)$ is visualized by an arrow drawn at the vertex 
corresponding to the component $E$ of  $\mu\inv (0)$ that $C$ intersects. 
% For each vertex $E$ the edges and arrows drawn from $E$ in the same anticlockw%ise order 
%as the corresponding curves cut $E$ (see the example below).   
To each vertex $E$ we assign its parity $p_\mu (E)$ and its multplicity 
$m (E)$.  
The multiplicity of 
the vertex is the generic multiplicity of $f\circ \mu$ on $E$.  
The parity of $E$ is $0$ if it has an orientable neighbourhood in 
$M$.  If a tubular neighbourhood of $E$ is a M\"obius band then its 
parity is $1$.

\begin{example}\label{examplesimple}
Let  $f(x,y) = x(x^3 - y^5)(x^3 + y^5)$ and 
$g(x,y) = x(x^3 - y^5)(x^3 - 2y^5)$.  The resolution graphs of $f$ and
$g$ are the following.   
\vspace{5mm}
\epsfxsize=14cm
\epsfysize=4cm
$$\epsfbox{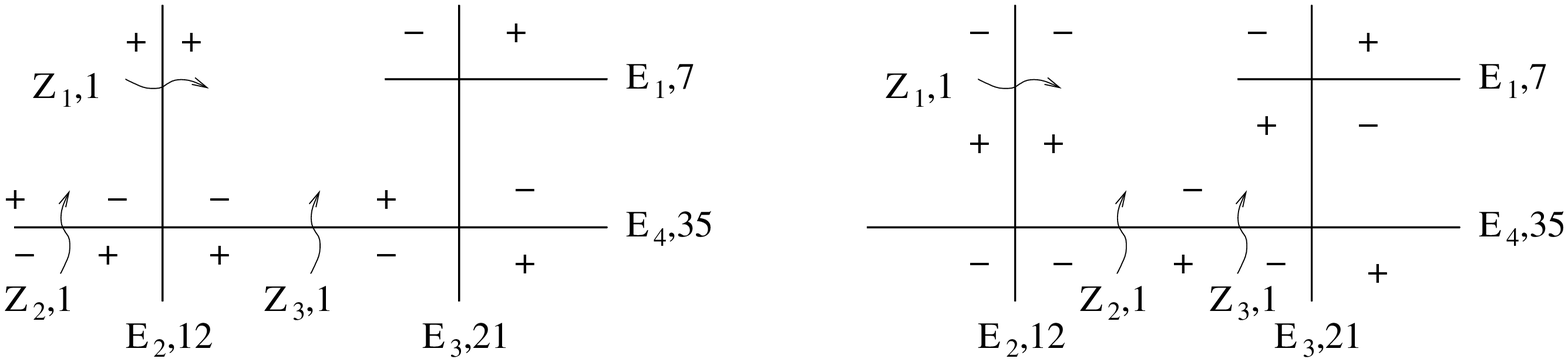}$$
\centerline{Resolution graph of $f$ \hspace{35mm} Resolution graph of
  $g$ \hspace{5mm}}\\

Note that $f$ and $g$ are desingularised by the same composition of
four point blowings-up $\mu$.  We denote by $E_i$'s the components of 
exceptional divisor of $\mu$, and by $Z_j$'s the components of the
strict 
transforms of $f\inv (0)$ and $g\inv (0)$ by $\mu$.  On the graph the numbers next
to these components denote their multiplicities.  The dual graphs of minimal resolution of $f$ and $g$
coincide.    
\vspace{5mm}
\epsfxsize=10cm
\epsfysize=2.5cm
$$\epsfbox{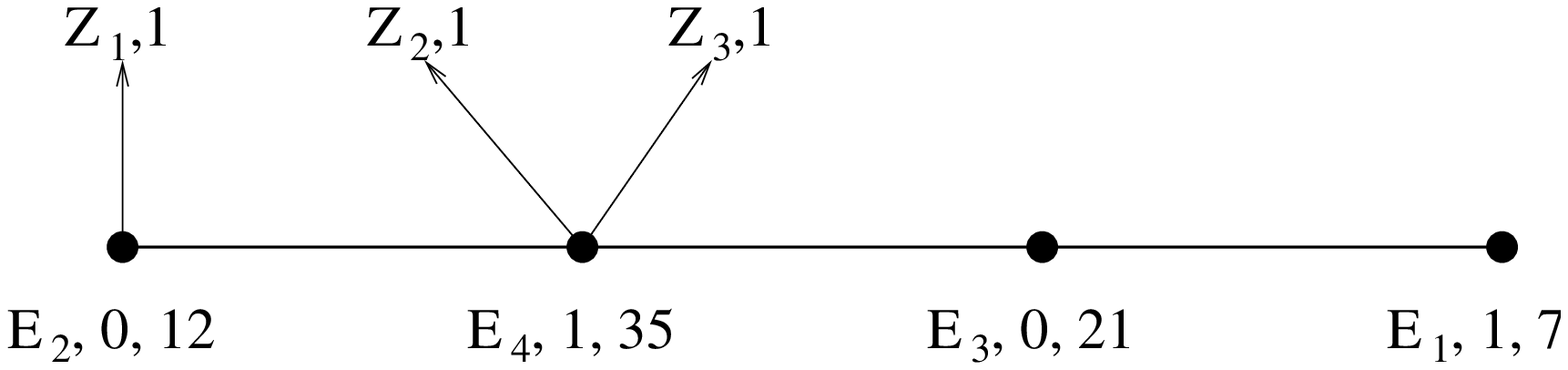}$$
%\centerline{$\R T_{(0,1)}(f)$\hspace{40mm} $\R T_{(0,1)}(g)$\hspace{5mm}}
%\vspace{5mm}
%\vspace{5mm}

The resolution graphs of $f$ and $g$ give  the same Fukui invariants and zeta functions,
cf. Theorem I \& VII in \cite{izumikoikekuo} and formulae (1.1) \&
(1.2) in \cite {koikeparusinski}. Therefore $f$ and $g$ have the same
Fukui invariants and zeta functions.  
As follows from Theorem \ref{allequivalent}, see also Example \ref{tree} below, 
 $f$ and $g$  are not blow-analytically equivalent.  
\end{example}

%%%%%%%%%%%%%%%%%%%%%%%%%%%%%%%%%%%%%%%%%%%%%%%%%%%%%%%%%%%%%%%

\subsection{Invariants of blow-analytic equivalence.}\label{invariants} 

We introduce a refinement of Fukui invariant that we will need later.  
First we recall briefly the construction Fukui invariant, cf.  \cite{fukui}.
Let $f : (\R^n,0) \to (\R,0)$ be an analytic function germ.
Set
$$
A(f) := \{ \ord (f(\gamma (t))) \in \N \cup \{ \infty \} 
; \, \gamma : (\R,0) \to (\R^n,0) \text { real analytic} \} .
$$ 
Let $\lambda : U \to \R^n$ be an analytic arc with $\lambda (0) = 0$, 
where $U$ denotes a neighbourhood of $0 \in \R$. 
We call $\lambda$  {\it nonnegative} (resp. {\it nonpositive}) 
{\it for} $f$ if 
$(f \circ \lambda)(t) \geq 0$ (resp. $\leq 0$) 
in a positive half neighbourhood $[0,\delta) \subset U$. 
Then we set
 
\vspace{3mm}

\qquad $A_+(f) := \ \{ \ord (f \circ \lambda) ; \, \lambda$
is a nonnegative arc through $0$ for $f \}$, 

\qquad $A_-(f) := \ \{ \ord (f \circ \lambda) ;  \, \lambda$ 
is a nonpositive arc through $0$ for $f \}$.

\vspace{3mm}

\noindent Fukui in \cite{fukui} proved that if analytic functions $f, g : (\R^n,0) \to (\R,0)$
are blow-analytically equivalent, then $A(f) = A(g)$,
$A_+(f) = A_+(g)$ and $A_-(f) = A_-(g)$.
We call $A(f)$, $A_{\pm}(f)$ {\em the Fukui invariant},
{\em the Fukui invariants with sign}, respectively.  The proof of 
Fukui was based on the fact that a blow-analytic homeomorphism sends 
a real analytic arc to a real analytic arc, that he proved only for 
these blow-analytic homeomorphism $h$ defined by the diagram \ref{baequivalence} 
with $\mu$ and $\tilde \mu$ blowings-up of coherent ideals.  We complete the general 
case in Corollary \ref{preservingarcs} below.  
Apart from the Fukui invariants, motivic type invariants
are introduced in \cite{koikeparusinski} and \cite{fichou}.

%We next introduce more refined blow-analytic invariants.
%We keep the same notations $f$ and $\lambda$ as above.
Let $C_+(f)$ (resp. $C_-(f)$) be the set of local connected components of
$\{ x \in \R^n ; f(x) > 0 \}$ (resp. $\{ x \in \R^n ; f(x) < 0 \}$)
as set-germs at $0 \in \R^n$. 
Let $C_+(f) = \{ V_1(f), \cdots , V_v(f) \}$ and
$C_-(f) = \{ W_1(f), \cdots , W_w(f) \}$.
For $i = 1, \cdots , v$ (resp. $j = 1, \cdots , w$),
we call $\lambda$ a {\it positive} (resp. {\it negative}) {\it arc in}
$V_i(f)$ (resp. $W_j(f)$) {\it for} $f$ if $(f \circ \lambda )(t) > 0$
(resp. $(f \circ \lambda )(t) < 0$) and $\lambda (t) \in V_i(f)$
(resp. $\lambda (t) \in W_j(f)$) over an open interval 
$(0, \delta ) \subset U$.
Then we set

\vspace{3mm}

\qquad $A_+^i(f) := \ \{ \ord (f \circ \lambda ) ;  \, \lambda$   
is a positive arc in $V_i(f)$ for $f \}$,

\qquad $A_-^j(f) := \ \{ \ord (f \circ \lambda ) ;  \, \lambda$
is a negative arc in $W_j(f)$ for $f \}$,

\vspace{3mm}

\noindent for $i = 1, \cdots , v$, $j = 1, \cdots , w$.
Using the argument of Fukui \cite{fukui} and Theorem \ref{rmodproperties} (c) below,
the collection of sets $A_+^i(f)$'s and $A_-^j(f)$'s is a blow-analytic invariant.
We call them {\it the refined Fukui invariants with sign}. {\it The
  refined Fukui invariant} (without sign) is defined similarly.

%%%%%%%%%%%%%%%%%%%%%%%%%%%%%%%%%%%%%%%%%%%%%%%%%%%%%%%%%%%%%%%

\subsection{Constructing blow-analytic equivalence.}\label{constructing}

\smallskip
\begin{defn}\label{weaklyisomorphic}
We say that two real analytic functions $f(x,y)$ and $g(x,y)$ have 
 \emph{weakly isomorphic resolution spaces} if there exist resolutions 
$\mu : M \to \R^2$, 
$\tilde \mu : \tilde M \to \R^2$ of $f(x,y)$ 
and $g(x,y)$, respectively,
and an analytic isomorphism $\Phi : M\to \tilde M$ 
such that:
\begin{enumerate}
\item
$\Phi$ conserves the exceptional sets and the strict transform 
of the zero sets: 
$$ \Phi  (\mu \inv (0)) = \tilde \mu \inv (0),\qquad  
\Phi ((f\circ \mu )\inv (0)) = (g\circ\tilde \mu )\inv (0)$$
\item  
$\Phi$ conserves the multiplicities: If $C$ is a component of  
$(f\circ \mu )\inv (0)$  then $\mult _C f\circ \mu = 
\mult _{\phi(C)} g\circ \tilde \mu $.
\item
 $\Phi$ conserves the signs: $f\circ \mu (p) > 0$ iff $g\circ \tilde \mu 
(\Phi (p)) > 0$.
\end{enumerate}
\end{defn}

The following result follows from Proposition 7.2 of \cite{fukui}.

\begin{prop}\label{fukuiprop}
If $f(x,y)$ and $g(x,y)$ have weakly isomorphic resolution spaces then 
$f$ and $g$ are blow-analytically equivalent. 
\end{prop}

\begin{rem}\label{Nash}
In Definition \ref{weaklyisomorphic} it is enough to assume that $\Phi$ is $C^\infty$ 
or even that it is only a homeomorphism. Indeed, a homeomorphism $\Phi$ satisfying 
(1)-(3) of Definition \ref{weaklyisomorphic} can be 
approximated by an analytic isomorphism with the same properties.  
This is fairly easy to see since we are in the two variable case.  
We can also use the following much deeper Nash approximation argument. 

Recall that Nash maps are real analytic maps with semi-algebraic graphs.   
Firstly, by a result of Shiota \cite{shiota} 
we may suppose that $f$ and $g$ are polynomial functions and hence 
that $M$ and $\tilde M$ are Nash manifolds, and that the components
 of the exceptional divisors and 
the strict transforms of $f\inv (0)$, $g\inv (0)$ are their Nash submanifolds 
resp., intersecting 
transversally.  Then the existence of $\Phi$ that is a Nash isomorphism and satisfies  
(1)-(3) of Definition \ref{weaklyisomorphic} follows directly from the proof of 
the Nash Isotopy Lemma of \cite{fukuikoikeshiota}.     

The above argument and Theorem \ref{allequivalent} show that polynomial function germs 
$f(x,y)$, $g(x,y)$ (and more generally Nash function germs), are blow-analytically 
equivalent if and only if 
they are blow-Nash equivalent, see also \cite{fichou}, \cite{fukuipaunescu2}.  
\end{rem}

%\smallskip
%%%%%%%%%%%%%%%%%%%%%%%%%%%%%%%%%%%%%%%%%%%%%%%%%%%%%%%%%%%%%%%
%%%%%%%%%%%%%%%%%%%%%%%%%%%%%%%%%%%%%%%%%%%%%%%%%%%%%%%%%%%%%%%%%%%%%%%%%%%

%%%%%%%%%%%%%%%%%%%%%%%%%%%%%%%%%%%%%%%%%%%%%%%%%%%%%%%%%%%%%%%%%%%%

\section{Real modifications}  
\label{realmodifications}
\smallskip

We recall this classical definition of real modification, cf. \cite{kuo8}, \cite{fukuikoikekuo}.  
We also introduce a slightly more general notion of singular real modification in order to 
have a notion that is stable by taking the strict transforms by blowings-up 
with smooth centre, see (b) of Theorem \ref{rmodproperties}.

\begin{defn}\label{realmodification}  Let $Y$ be a real 
analytic manifold of pure dimension $n$.  We say that $\sigma:X\to Y$ 
is  a {\em singular real modification} if the following property is
satisfied. 

$X$ is a real analytic space,  $\sigma:X\to Y$ is   
 a proper surjective real analytic map, and there 
exist complexifications $\compl X, \compl Y$ of $X$ and $Y$, 
respectively, and a holomorphic extension 
$\compl {\sigma} : \compl X \to \compl Y$ of $\sigma$, that satisfy: 

$\compl X$ is a complex analytic space 
of pure complex dimension $n$ and  
 $\compl {\sigma}$ is an isomorphism in the complement of a  
closed nowhere dense subset $B$ of $\compl X$. (that is 
$\compl {\sigma}$ resticted to $\compl X\setminus B$ is open and 
an isomorphism onto its image.)
 
If, moreover, $X$ is nonsingular, then 
 we say that $\sigma$ is a {\em real modification}.   
\end{defn}

 Note that a real analytic map that 
is an isomorphism in the complement of a closed nowhere dense subset of $X$ is 
not necessarily a real modification.  This is for instance the case for
 $\sigma:\R\to \R$ given by $\sigma(x)=x^3$.  

\begin{rem}
If one defined  real modifications  simply 
 as compositions of blowings-up with smooth centres,  then one 
would need the strong 
factorisation in order to show that the induced notion of blow-analytic 
equivalence, see subsection \ref{b.a.e.},  is an equivalence relation.   
At the moment we do not know whether the strong factorisation holds. 
With Kuo's definition, certainly quite complicated and less natural, the proof of transitivity of the blow-analytic equivalence is fairly easy.  
\end{rem}

\begin{thm}\label{rmodproperties}
Let $\sigma:X\to Y$ be a singular real modification.  Then 
\begin{enumerate}
\item[{ }(a)]
$\sigma$ is an isomorphism 
over the complement in $Y$ of a subanalytic subset $A \subset Y$ 
of real codimension $2$.  
\item[{ } (b)] 
%The family of real singular modifications is stable by taking the 
%strict transform by local blowings up:  
Let $\pi :Y'\to Y$ be a blowing-up with a nonsingular nowhere dense centre.    
Then the strict transform 
$\sigma':X'\to Y'$ by $\pi$ 
is a singular real modification.   
\item[{ }(c)]
If $\gamma : (\R,0) \to (Y,p)$ is the germ of a real analytic arc at 
$p\in Y$ then there is a real analytic  
$\tilde \gamma : (\R,0)  \to (X,\tilde p)$, $\tilde p\in X$,  
such that 
$\sigma \circ \tilde \gamma = \gamma$.  \\
Moreover, there is a closed subanalytic nowhere dense 
$A\subset Y$ (independent of $\gamma$) such that if $\gamma \inv (A)$ is discrete 
then such $\tilde \gamma$ is unique. \\
(We say for short that $\sigma$ satisfies \emph {the arc 
lifting property} and   \emph {the unique generic arc 
lifting property} ). 
\end{enumerate}
\end{thm}

\begin{proof}
Let $\sigma :X\to Y$ be a singular real modification and let $B$ 
be a closed nowhere dense subset of $X_{\C}$ sucht that $\sigma _{\C}$ is an 
isomorphism in the complement of $B$.  

First we show (b). Let 
$\pi: Y' \to Y$ be a 
blowing-up with smooth nowhere dense centre $C$ and consider the diagram 
\begin{equation}\label{stricttransform}
\minCDarrowwidth 1pt
{\begin{CD}
 X @< \pi' << X' \\ 
@V\sigma VV @VV{\sigma}'V \\
Y @<\pi<< Y'  
\end{CD}}
\end{equation}  
where $\sigma':X'\to Y'$ is the strict transform of $\sigma$ by $\pi$.  
Then $\pi'$ is the blowing-up of the pullback 
by $\sigma$ of the ideal of $C$.   The above diagram induces 
the complexified diagram, where we complexify $\pi$ and $\pi'$ 
by the corresponding complex blowings-up. Let $E'_\C$ denote the exceptional 
divisor of $\pi'_\C$.   Then $\compl {\sigma'}$ 
is an isomorphism in the complement of $E'_\C \cup \compl{ \pi'} \inv (B)$ 
that is nowhere dense in $X'_{\C}$.    

To show (a) we note first that by assumption the generic fibres of 
$\sigma$ consist of single points.  
Over  a generic point of $y\in Y$ in codimension $1$, $\sigma_\C$ is finite. 
Indeed,
let $y= \sigma(x)$ and suppose that $(\sigma_\C)_x$ were not finite.  
Then $X _\C$ 
would contain a vertical component $(X_\C)_1$   
(i.e. on such  component  
the rank of the differential of $\sigma_\C$ restricted to
 the regular part of 
$(X_\C)_1$ is everywhere smaller than $n$). 
But this contradicts the existence of a nowhere dense 
$B$ of the definition of real modification.  Therefore over 
$y$, $\sigma_\C$ is finite.  
But a finite real modification has to be  of degree one, that is an isomorphism.  

If $\sigma$ is a blowing-up with a smooth nowhere dense centre, or a composition 
of such blowings-up, then it satisfies (c).  Thus (c) for an arbitrary real modification 
follows from the local flattening theorem, cf. \cite{hironaka1}, 
\cite {hironakalejeuneteissier}.  Indeed, for any germ $\gamma : (\R,0) \to (Y,p)$ there 
is a composition of local blowings-up 
with smooth centres $\pi: Y' \to Y$ such that $\gamma$ lifts to $Y'$ and that
 the strict transform $\sigma': X'\to X$  of $\sigma$ by $\pi$ is flat.  Since $\sigma'$ is 
a singular real modification by (b) it has to be an isomorphism.  Therefore $\gamma$ lifts 
to $X'$ and hence to $X$, see \eqref{stricttransform}.  
The last claim of (c) follows now easily from (a). 
\end{proof}

\begin{rem}
The arc-lifting property for real modifications has been proven also in 
\cite{fukuipaunescu2}, section 5. 
\end{rem}

\begin{cor}\label{preservingarcs}
Let $h:(\R^n,0) \to (\R^n,0)$ be a blow-analytic homeomorphism and let 
$\gamma : (\R,0) \to (\R^n,0)$.  
Then $\gamma$ is real analytic if and only if so is $h\circ \gamma$.   
\end{cor}

\smallskip
\begin{proof}[Proof of Theorem \ref{2variablemod}.] 

Clearly a composition of point blowings-up is a real modification.  We shall 
show the converse.  Suppose that $\sigma$ is a real modification. 
By  (a) of 
Theorem \ref{rmodproperties} there exists a discrete subset  
$A\subset S$ such that $\sigma$ is an isomorphism 
over $S\setminus A$.  Fix $p\in A$ and suppose that $\sigma$ 
is not an isomorphism over $p$. Let $\pi  :\hat S\to S$ be a blowing-up of 
$S$ at $p$.  
The main point is to show that $\sigma$ factors through $\pi $, 
that is that there is $\hat \sigma:X \to \hat S$ 
such that $\pi \circ \hat \sigma =\sigma$.  The proof follows 
a classical argument of  
elimination of indeterminacy of  rational maps between algebraic surfaces, 
see e.g. 
\cite {beauville} Proposition II.8.  
Since the problem is local, we shall work in a neighbourhood of $p$ 
and assume that $\sigma$ 
is an isomorphism over $S\setminus \{p\}$.    

\begin {lem}\label{zariski}
 Let 
$q\in \sigma\inv (p)$ and assume that $\sigma$ is not an isomorphism at $q$.  
Let $x,y$ be a system of  
 local analytic coordinates $x,y$ at $p$
Then there is a linear combination $z=ax+by \in m_p \setminus m_p^2$  such
that $\sigma^*(z) \in m_q^2$.  
\end{lem}

\begin {proof}
Indeed, by assumption, $\sigma$ is not an isomorphism at $q$ so $\sigma^*(x)$  
and $\sigma^*(y)$ are linearly dependent in $m_q/m_q^2$.
\end{proof}

Let $\hat \pi:\hat X \to X$ be the blowing-up of the ideal $\sigma^*(m_p)$.  
%and let $\hat \sigma:\hat X \to \hat S$ be the induced map.  
Then $\hat X$ can be identified with this irreducible component of 
$X\times _S \hat S$ that is not entirely included in the inverse image 
of the exceptional divisor $E$ of $\pi$. 
%, and  $\hat \sigma$ is induced by the projection on the second factor.  
\begin{equation*}\minCDarrowwidth 1pt
\begin{CD}
 X @< \hat \pi << \hat X @.\subset X\times _S \hat S \\ 
@V\sigma VV @VV\hat \sigma V @.\\
S @<\pi<< \hat S @. 
\end{CD}
\end{equation*}  
 Note that 
   $\hat \pi:\hat X \to X$ is surjective since the zero set of
   $\sigma^*(m_p)$ is nowhere dense in $X$.  
By (a) of Theorem \ref{rmodproperties}, $\hat  \sigma $ is 
an isomorphism over the complement of a finite 
subset $F$ of the exceptional divisor $E\subset \hat S$ and hence 
the map $\hat  \sigma \inv :\hat S \setminus F\to \hat X$ 
is well-defined.   

\begin{lem}
$\hat \pi \circ \hat  \sigma \inv :E \setminus F\to  \sigma \inv (p)$ 
is not constant.   
\end{lem}

\begin{proof}
Suppose that this map is constant and that its image is $q \in \sigma \inv (p)$.   
For any $\hat p\in E\setminus F$ and for $z\in m_p$ given by lemma 
\ref{zariski}
$$
(\hat \pi \circ \hat \sigma\inv )^* \sigma^* (z) \in m_{\hat p}^2.
$$
But this is impossible since 
$$
(\hat \pi \circ \hat \sigma\inv )^* \sigma^*(z) 
= \pi^*(z) \notin  m_{\hat p}^2
$$
for all $\hat p \in E$ but one (the one corresponding to the zero 
set of $z$).  
\end{proof}

Thus for any $q\in \sigma\inv (p)$, $\hat \pi \inv (q)$ is finite.  
It follows from the next lemma that $\hat \pi$ has to be an 
analytic isomorphism.

\begin{lem}
Let $X$ be a nonsingular real analytic surface and let 
$\Pi: \mathcal X \to X$ be the blowing-up of an ideal $\mathcal I$ 
(not identically equal to zero).  
Suppose that locally at any point of $X$, $\mathcal I$ 
can be  generated by two real analytic functions.   
Then, if all the fibres of $\Pi$ are finite  then 
$\Pi$ is an analytic isomorphism.
\end{lem}

\begin{proof}
We work locally on $X$ so we assume that $\Pi$ is an isomorphism over the 
complement of a single point $q\in X$.  
Since $\mathcal I_q$ has two generators, 
$\mathcal X \subset X\times \projr^1$.  
Then, because of dimensional reason, the geometric finiteness of $\Pi$
implies that its complexification is also finite.  More precisely, 
fix a $\hat q \in \Pi \inv (q)$ 
and work locally in a neighbourhood of $\hat q$.    
The complexification $(\compl {\mathcal X})_{\hat q}$ of 
$\mathcal X_{\hat q}$  
is a complex analytic subset of $(\compl {X}\times \projc^1)_{\hat q}$, 
where by 
$\compl {X}$ we denote the complexification of $X$ at $q$.  
Denote by $\compl {\Pi}:(\compl {\mathcal X})_{\hat q}\to (\compl {X})_q$
 the projection onto the first factor.  Then 
$\compl {\Pi}\inv (q)$ is a complex analytic subset of 
$(\{ q \} \times \projc^1)_{\hat q}$ 
whose real part is reduced to one 
point and hence is finite itself.  Thus  $(\compl {\Pi})_{\hat q}$ is a finite 
map.  Since a finite blowing-up of a nonsingular complex analytic
space is an isomorphism, so is $\Pi$.  
\end{proof}

Thus we have shown that $\sigma$ factors through $\pi$ as claimed.  
Then we apply the same procedure to $\hat \sigma$.  
To finish the proof, we note  that 
after a finite number of point blowings-up (more precisely locally 
finite on $S$) this process terminates, i.e. the obtained lift 
is an isomorphism.  Indeed, after each blowing-up the number of 
irreducible components of the exceptional set  increases but it 
 cannot be bigger than the number of 
irreducible (global, analytic) ones of $\sigma \inv (p)$.  The  proof
of 
Theorem \ref{2variablemod} is complete.  
\end{proof}

%%%%%%%%%%%%%%%%%%%%%%%%%%%%%%%%%%%%%%%%%%%%%%%%%%%%%%%%%%%%%%%%%%%

%%%%%%%%%%%%%%%%%%%%%%%%%%%%%%%%%%%%%%%%%%%%%%%%%%%%%%%%%%%%%%%%%%%%

\section{Proof of (1)$\Rightarrow$(2) of 
   Theorem \ref{allequivalent} }
\label{lift}
\smallskip

Let $f : (\R^2,0) \to (\R,0)$ and $g : (\R^2,0) \to (\R,0)$ be 
real analytic function germs.  By Theorem  \ref{2variablemod} there
exists a commutative diagram
\begin{equation}\label{exists?h1} 
\minCDarrowwidth 1pt
{\begin{CD}
(M, {\mu'} \inv (E)) @>\mu'>> (\wR ,E_1) @>\pi>> (\R^2,0) @>f>>\R  \\ 
@V\Phi VV @V\exists  \, h_1 \, ?VV @ VhVV    @| \\
(\tilde M, \tilde {\mu'} \inv (E)) @>\tilde {\mu}'>> (\wR , \tilde E_1) 
@> \tilde\pi>> (\R^2,0) @>g>>\R   
\end{CD}}
\end{equation} 
where $\pi$ and $\tilde \pi$ is the 
blowing-up of the origin, $\mu'$ and $\tilde \mu'$ are compositions of
point blowings-up, 
$\Phi$ is 
an analytic isomorphism and $h$ is a homeomorphism such that $f = g \circ h$.

\begin{prop}\label{h1exists}
There exists a homeomorphism 
$h_1: (\wR ,E_1)  \to (\wR,\tilde E_1) $, 
$\tilde \pi \circ h_1 = h\circ \pi$, 
closing the diagram \ref{exists?h1}.   
\end{prop} 

\begin{proof}
(For simplicity we shall drop the germ notation.)   
By composing $\mu = \pi \circ \mu'$ and $\tilde \mu = \tilde \pi \circ \tilde \mu'$ with 
additional blowings-up, if necessary, we may assume that they are resolutions 
of $f$ and $g$ respectively, that is 
$f\circ \mu$ and $g\circ \tilde \mu$ have only normal crossing singularities.   

Denote by $E_1, \ldots, E_k$ the components of the exceptional 
divisor of $\mu$, in the order they were created 
(we shall keep the same notation for the exceptional divisor of a point blowing-up 
and for its subsequent strict transforms). We show below that 
 among these components, $E_1$ can be recognised on the dual
 resolution graph
 $S(\mu )$
 of $\mu$.  
Hence we  conclude that $\Phi (E_1)=\tilde E_1$, that is the most difficult step 
in showing the existence of $h_1$.

\begin{lem}\label{charE1}
Let $m=\mult_0 f$ and suppose that the zero set of the leading homogeneous part $f_m$ 
of $f$ is not reduced to a point.  Then among all divisors of multiplicity $m$ in 
$\mu \inv (0)$, $E_1$ is completely characterised by the following property: 
\begin{enumerate}
%\item [] 
\item[{ }(H)]
Either the strict transform of $f\inv (0)$ intersects $E_1$ 
or there is a connected component of 
$S(\mu) \setminus E_1$ that does not contain a vertex of multiplicity $m$.
\end{enumerate}
(The two conditions of (H) do not exclude themselves and 
often they both hold for $E_1$.)
\end{lem}

\begin{proof}
We first show that property (H) is hereditary for the divisors of multiplicity $m$ 
in the following sense.  Suppose that $\sigma :N \to \R^2$ is a real modification 
and let $E$ be an exceptional 
divisor of $\sigma$.  Suppose moreover that 
$f\circ \sigma$ is a normal crossings 
in a neighbourhood of $E$.  Let $\sigma':N'\to N$ be a blowing-up of $p\in E$.  
Then we say that (H) is \emph{hereditary} if $E$ satisfies (H) as a divisor of $N$ 
iff its strict transform satisfies (H) as a divisor of $N'$.  
In our case one may check the heredity easily by inspection, since  either 
$\mult_p f\circ \sigma =m$ or $E$ intersects at $p$ another divisor or 
the strict tranform of $f\inv (0)$ and then $\mult_p f\circ \sigma > m$ .    

In general, let $E$ be a component of the exceptional divisor
 of a real modification $\sigma :N \to \R^2$.  
We call a point $p\in E$ \emph{simple} (with respect to $f$) if 
$\mult_p f\circ \sigma =\mult_E(f)$.  
A divisor $E\ne E_1$ of multiplicity $m$ can be  created only by blowing-up 
a simple point on another divisor of multiplicity $m$. Therefore for such 
a divisor $f\circ \sigma $ is always normal crossings in a neighbourhood of $E$ and 
does not satisfy (H) at the moment it is created,  and hence, by heredity, never.     

By assumption on $f_m$,  the divisor $E_1$, when created by the first blowing-up $\pi$, contains 
 some points of multiplicity higher than $m$. 
If its strict transform in $M$ does not intersect the strict transform of 
$f\inv (0)$, it means that all the points of multiplicity higher than $m$ on $\pi \inv (0)$ 
have been blown-up.  
Each of such blowing-up produces a divisor of multiplicity higher 
than $m$ and a connected 
component of $S(\mu) \setminus E_1$ that does not contain 
a divisor of multiplicity $m$. 
Blowing-up points on the divisors on this component 
cannot produce new divisors  of multiplicity $m$.  
This ends the proof of lemma.   
\end{proof}
 
Suppose that $f_m\inv(0) = \{0\}$. Then the blowing-up of the origin 
$\pi:\wR \to \R^2$ resolves $f$. 
The modification $\mu$ is the composition of $\pi$ and finitely many 
point blowings-up.  To each component $E$ of  $\mu\inv (0)$ we
associate the following number 
\begin{equation}\label{graphinvariant}
\varepsilon_\mu (E) = \sum_{E' \in \mathcal E_\mu} \frac {m(E')}{m(E)} 
(E',E)_\mu,   
\end{equation}
where by $\mathcal E_\mu$ we denote the set of all exceptional divisors of 
$\mu$ and $(E,E')_\mu$ is defined by the following rule 
\begin{equation}\label{intnumber}
(E,E')_\mu = 
\begin{cases}
1 \qquad \text { if } E\ne E'\text { and } E\cap E' \ne \emptyset \\
0 \qquad \text { if } E\cap E' = \emptyset \\
p_\mu(E) \quad \text { if } E=E' ,  
\end{cases}
\end{equation}
where the parity $p_\mu(E)$ of $E$ equals $0$ if $E$
admits an orientable neighbourhood in $M$, otherwise $p_\mu(E)=1$.  
That is $(E,E)_\mu \pmod 2$ equals the intersection number 
of  $E$ and $E'$. 
  We can detect $E_1$ among all the other components thanks to the 
following 
lemma. 

\begin{lem}\label{charE2}
Under the above assumptions $\varepsilon_\mu (E)\in \Z$ and 
$$
\varepsilon_\mu (E_1) \equiv 1 \mod 2, 
\qquad \varepsilon_\mu (E) \equiv 0 \mod 2 
\text { for } E\ne E_1.
$$
\end{lem}

\begin{proof}
We check how $\varepsilon_\mu (E)$ changes under blowings-up.
Let $\sigma:N\to M$ be a blowing-up of $p\in E$ and 
$\mu' = \mu \circ \sigma$.  
If $p$ is a simple point  of $E$ then a new divisor of multiplicity $m(E)$ 
is produced and the parity of $E$ changes.  Both these events affect the sum in 
\eqref{graphinvariant} by adding $\pm 1$.  Similarly, if $p\in E\cap E'$ then 
a new divisor $E"$ of multiplicity $m(E) +m(E')$ is created.  Then 
$(E,E")_{\mu'}=1$ and  $(E,E')_{\mu'}=0$.  
Moreover, the parity of $E$ is reversed.  Therefore  
$$
\frac {m(E")}{m(E)} (E,E")_{\mu'} + \frac {m(E)}{m(E)} 
(E,E)_{\mu'} \equiv  \frac {m(E')}{m(E)} (E,E')_{\mu} + \frac {m(E)}{m(E)} 
(E,E)_{\mu} \mod 2
$$
that shows that $\varepsilon_\mu (E)  \mod 2$ does not depend on $\mu$.
Then we compute $\varepsilon_\mu (E)  \mod 2$ at the moment $E$ is created.  
If $E=E_1$ we take $\mu = \pi$ and get 
$$\varepsilon_\pi (E_1) = 1,$$
as claimed.  
If $E\ne E_1$ we denote by $\sigma :N \to \R^2$ the real modification that has 
$E$ as the last created divisor.  Then $E$ is the exceptional divisor of the 
blowing-up of either a simple point on a divisor $E_i$,,  or of the 
intersection point of 
two distinct divisors $ E_i,  E_j$.  In the former case 
$$\varepsilon_\sigma (E) = \frac {m(E_i)}{m(E)} 
(E_i,E)_{\sigma }  + \frac {m(E)}{m(E)} 
(E,E)_{\sigma} = 1 + 1 \equiv 0 \mod 2.  $$
In the latter case 
\begin{eqnarray*}
& \varepsilon_\sigma (E) & = \frac {m(E_i)}{m(E)} (E_i,E)_{\sigma } + 
\frac {m(E_j)}{m(E)} (E_j,E)_{\sigma }  
+ \frac {m(E)}{m(E)} (E,E)_{\sigma}  \\
& & = \frac {m(E_i)}{m(E)} + \frac {m(E_j)}{m(E)} + 1  \equiv 0 \mod 2.  
\end{eqnarray*} 
This ends the proof of lemma.  
\end{proof}

We continue the proof of Proposition \ref{h1exists}.  If $f_m\inv (0) = \{0\}$, 
then $E_1$ does not satisfy the property (H) of Lemma \ref{charE1}.  
Therefore there is no exceptional divisor of $\mu$ of multiplicity $m$ that 
satisfies this property.  Consequently the same is true for $\tilde E_1$ and 
$\tilde \mu$.  Therefore $g_m\inv (0) = \{0\}$ and we may conclude by Lemma 
\ref{charE2} that $\Phi (E_1)=\tilde E_1$.  

Similarly if $E_1$ satisfies the property (H) then so does $\tilde E_1$ and 
$\Phi (E_1)=\tilde E_1$ by Lemma 
\ref{charE1}. 

Denote by $F\subset \pi \inv (0)$, resp. $\tilde F\subset \tilde \pi \inv (0)$,  
the image of the other  exceptional divisors of $\mu$, resp. 
 $\tilde \mu$.  Thus $F$ and $F_1$ are finite. 
For each $p\in F$, $(\mu') \inv (p)$ is the union of all divisors in a connected 
component of $S(\mu) \setminus E_1$.  Therefore $\Phi$ sends $(\mu') \inv (p)$  onto 
$(\tilde {\mu}') \inv (q)$ for a unique $q\in \tilde F$.  
Moreover $\Phi$ induces a homeomorphism of 
$E_1\setminus F$ onto $\tilde E_1\setminus \tilde F$.   Since both $\mu'$ and 
$\tilde \mu'$ are proper this is sufficent to  conclude that $\Phi$ 
induces a homeomorphism $h_1: (\wR ,E_1)  \to (\wR,\tilde E_1) $, as claimed.  
This ends the proof of Proposition \ref{h1exists}.  
\end{proof}

\begin{lem}\label{normalcrossing}
If $f : (\R^2,0) \to (\R,0)$ is normal crossing and 
$f$ and $g : (\R^2,0) \to (\R,0)$ are 
blow-analytically equivalent then $g$ is also normal crossing 
(with the same exponents).
\end{lem}

\begin{proof}
We shall only consider the case $f=x^ay^b$, $a>0,b>0$, 
the case when $f=x^a$ is 
easier.  Then the zero set of $g$ contains two real analytic curves, one with generic 
multiplicity $a$ and the other one with generic multiplicity $b$. Therefore 
$$
g = (g_1)^a (g_2) ^b g_3.
$$
Since the multiplicity is a blow-analytic invariant, $\mult_0 g = a+b$. 
 Consequently $g_3$ is a unit and $g_1,g_2$ are regular germs.  

The zero set $f\inv (0)$ of $f$ divides a neighbourhood of $0$ in $\R^2$ into 
4 sectors and each of them contains a half branch of real analytic curve 
$\gamma (t): (\R, 0)\to (\R^2, 0)$ on which $f\circ \gamma (t)$ has
order $a+b$.  The invariance of the refined Fukui invariant, cf. subsection
\ref{invariants}, shows that  the same is true for 
the four sectors of the complement of 
$g\inv (0)$.  This implies that the zero sets $g_1\inv (0)$  and 
$g_2\inv (0)$ are transverse and ends the proof of lemma.   
\end{proof}

\smallskip
Suppose  that the blow-analytic equivalence of $f$ and $g$ be given \eqref{baequivalence}. 
Performing additional point blowings-up, if necessary, we may assume that 
both $f\circ \mu$ and $g\circ \tilde \mu$ are normal crossings.  
By Theorem \ref{2variablemod}, 
$\mu$ is the compostion of a sequence of points blowings-up  
\begin{equation}%\label{exists?h1}
\minCDarrowwidth 1pt
\begin{CD}
M=M_k @>b_{k}>> M_{k-1}@>b_{k-1}>> \cdots @>b_{2}>> 
M_1 =\wR @>b_{1}=\pi >> \R^2 
\end{CD}
\end{equation}
The order of choice of centres of the blowings-up is not unique,  
if we have to blow-up two different points on $\wR$ we may do it in any order.  
Suppose that we first blow-up only singular points (in the sense of 
non-normal 
crossing) of the total transforms of 
$f$.  The composition of these, say $s(f)$, 
blowings-up is the minimal resolution of $f$.  
 Then, using repeatedly Proposition \ref{h1exists} 
we obtain  
\begin{equation}\label{cascade3}
\minCDarrowwidth 1pt
\begin{CD}
M @>b_{k}>> M_{k-1}@>b_{k-1}>> \cdots @>b_{2}>> M_1  @>b_{1} >> \R^2 @>f>> \R\\
@V\Phi VV  @ VVh_{k-1}V @. @ VVh_{1}V @ VVhV @ | \\
\tilde M @>\tilde b_{k}>> \tilde M_{k-1}@>\tilde b_{k-1}>> 
\cdots @>\tilde b_{2}>> \tilde M_1 
 @>\tilde b_{1} >> 
\R^2 @>g>> \R .
\end{CD}
\end{equation}
By Lemma \ref{normalcrossing}, 
whenever $p\in M_i$ is singular (in the sense of 
non-normal crossing) for  the total 
transform of $f$ then $h_i(p)\in \tilde M_i$ is also singular 
for the total transform of $g$, 
and vice versa.  Consequently, the composition of the first $s(f)$ 
blowings-up of the lower row of \eqref{cascade3} 
is the minimal resolution of $g$. 

This ends the proof of  (1)$\Rightarrow $(2) of 
Theorem \ref {allequivalent}.

%%%%%%%%%%%%%%%%%%%%%%%%%%%%%%%%%%%%%%%%%%%%%%%%%%%%%%%%%%%%%%%%%%%%%%%%

%%%%%%%%%%%%%%%%%%%%%%%%%%%%%%%%%%%%%%%%%%%%%%%%%%%%%%%%%%%%%%%%%%%

\medskip
\section{Real analytic demi-branches}
\label{arcs}
\smallskip

By {\it a (parametrised) real analytic  arc at $0\in \R^2$} 
we mean an analytic non-constant germ 
$\gamma (t) : (\R,0) \to (\R^2,0)$.  We consider two such arcs 
$\gamma, \tilde \gamma$ 
\emph {equivalent} if there is an analytic orientation preserving 
isomorphism germ $\sigma :  (\R,0) \to (\R,0)$ 
such that $\tilde \gamma = \gamma \circ \sigma$.  
The equivalence classes will be called 
\emph{real analytic arcs at $0\in \R^2$}.  
  
A real analytic arc $\gamma$ will be called \emph{reduced} 
if it cannot be written down as 
$\gamma (t) = \eta (t^k)$, where $k>1$ and $\eta$ is a real analytic arc.  
If $\gamma$ is reduced then it is injective as a map and its image is an 
 irreducible  analytic 
 germ of dimension $1$.  

\begin{rem}
If $h:(\R^2,0) \to (\R^2,0)$ is a blow-analytic homeomorphism and 
$\gamma (t) : (\R,0) \to (\R^2,0)$ is reduced then so is $h\circ \gamma$.
\end{rem}

By {\it a real analytic demi-branch at $0\in \R^2$} we mean a 
parametrised real analytic  
 arc $\gamma (t) : (\R,0) \to (\R^2,0)$ restricted to $t\ge 0$.  
We again identify   $\gamma$ and  
$\tilde \gamma$ if $\tilde \gamma = \gamma \circ \sigma$ 
for an orientation preserving analytic isomorphism 
$\sigma :  (\R,0) \to (\R,0)$.  
By Puiseux Theorem each  reduced real analytic demi-branch can be expressed,
 after a coordinate change in $(\R^2,0)$, as $\gamma(t) = (\lambda
 (t^m), t^m)$, 
where $\lambda$ is a fractional power series such that $\ord_0 \lambda (y)
\ge 1$.   If this is the case then we say, for short, 
that the demi-branch $\gamma$ is given by 
\begin{eqnarray}\label{fps}
x= \lambda (y) = a_{m_1 / m}  y^ {m_1 / m} + a_{m_2 / m}  y^{m_2 / m}
+ \cdots ,  \qquad y\ge 0
\end{eqnarray}
where $m \le m_1 < m_2 < \cdots$ are positive integers, having no
common divisor.  Sometimes we also assume that $m_1>m$ and that $m_1$ does
not divide $m$. For a single arc this is always possible after another 
analytic coordinate change.

By \emph {the tangent direction at $0$ of a demi-branch $\gamma (t)$} we mean 
$\lim _{t\to 0} \frac {\gamma(t)} {\|\gamma (t)\|}\in \rm S^1$.  
Let $\gamma_1, \gamma_2$ be reduced real analytic demi-branches 
tangent at $0$.   Suppose that $\gamma_1, \gamma_2$ are given, 
in the same system of coordinates by Puiseux series $\lambda_1$, 
resp. $\lambda_2$.  
We define \emph{the contact order of}  $\gamma_1$ and $\gamma_2$ as 
$$
O(\gamma_1,\gamma_2) := \ord_0 \, (\lambda_1 - \lambda_2)(y).
$$
If $\gamma_1, \gamma_2$  have distinct tangent directions at $0$ then we define 
$O(\gamma_1,\gamma_2) =1$.  The order contact 
between two real analytic demi-branches is well-defined, that is it is 
 independent 
of their Puiseux presentations of the form \eqref{fps}.

\smallskip
The Puiseux pairs of a reduced real analytic demi-branch $\gamma$ 
are pairs of relatively prime positive integers 
$(n_1,d_1), \ldots, (n_q,d_q)$, $d_i>1$ for  $i=1,\ldots,q$,   
$\frac {n_1}{d_1} < \frac {n_2}{d_1 d_2} < \cdots < 
\frac {n_q }{d_1\ldots d_q}$, 
such that  
\begin{eqnarray}\label{Ppairs} 
& & \qquad  \qquad  \qquad   \lambda (y)  =   \sum_\alpha a_\alpha y^\alpha = \\ 
\notag
&  & = \sum_{j=1}^{[{n_1}/{d_1}] }  a_j y^j  %  +a_{ {n_1}/{d_1}}y^{ {n_1}/{d_1}} 
+ \sum_{j=n_1}^{[{n_2}/{d_2}] }   a_{{j}/{d_1}}y^{ {j}/{d_1}} 
% + a_{ {n_2}/{d_1 d_2}}y^{ \frac {n_2} {d_1 d_2}} 
+ \sum_{j=n_2}^{[{n_3}/{d_3}] } a_{ {j}/{d_1 d_2}}y^{ \frac {j} {d_1 d_2}} + \cdots 
+ \sum_{j=n_q}^\infty   a_{{j}/{d_1 d_2 \cdots d_q}}y^{ \frac {j} {d_1 d_2 \cdots d_q }} 
\end{eqnarray}
and $a_{n_i/d_1\cdots d_i}\ne 0$ for $i=1,\ldots ,q$.  Following \cite{wall} 
chapter 3 we call the integers 
$m; \beta_1 =  n_1 d_2 \cdots d_q, \beta_2 = n_2 d_3 \cdots d_q, \ldots, \beta_q=n_q$ 
\emph{the Puiseux characteristic sequence of $\gamma$}. We will call
$m =d_1 d_2 \cdots d_q $ \emph{the multiplicity of $\gamma$}.  
The coefficients 
$A_i (\gamma ):= a_{n_i/d_1\cdots d_i}= a_{\beta_i/m}$ for $i=1,\ldots ,q$ 
will be called 
\emph{the characteristic coefficients of $\gamma$}.  

\begin{prop}\label{cs}
Let $\gamma$ be a reduced real analytic demi-branch of multiplicity $m$.  
  Then the signs of  characteristic coefficients are well-defined,
  they are independent 
of the Puiseux presentation  $\gamma(t) = (\lambda (t^m), t^m)$ 
in an oriented system of 
coordinates at $(\R^2,0)$.  \
\end{prop}

\begin{proof}
Let $\gamma(t) = (\lambda (t^m), t^m)$. 
Fix $i =1,\ldots , q$ and consider a new demi-branch 
$\gamma_i$ defined in the same system of coordinates by 
$$
\gamma_i (t) = (\lambda_i (t^m), t^m),\qquad 
\lambda_i  (y)  =   \sum_{\alpha<\beta_i/m} a_\alpha y^\alpha .
$$
The order of contact between 
$\gamma$ and $\gamma_i$ equals exactly $\beta_i/m$ and 
$(n_i,d_i)$ is not the i-th Puiseux pair of $\gamma_i$.  These two properties 
characterise $\gamma_i$  "up to terms of order higher than $\beta_i/m$''.  
Since $\beta_i/m > 1$ both demi-branches are tangent.  
The i-th characteristic coefficients of $\gamma$ is positive 
if and only if $\gamma$ follows $\gamma_i$ in the clock-wise direction in 
the orientation of $(\R^2,0)$ induced by the coordinate system $x,y$.   
%This ends the proof
\end{proof}

\smallskip
%%%%%%%%%%%%%%%%%%%%%%%%%%%%%%%%%%%%%%%%%%%%%%%%%%%%%%%%%%%%%%%
\subsection{Effect of a blowing-up.}

Let $\gamma (t)$ be a reduced real analytic demi-branch.  Denote by 
$\tilde \gamma$ the strict transform of $\gamma$ by the 
blowing-up of the origin $\pi: \wR  \to \R^2$.  Then $\tilde \gamma$ 
is a reduced real analytic demi-branch based at a point of $\wR$ that 
we denote by $\tilde 0$.  The  Puiseux characteristic sequence 
of $\tilde \gamma$ is the following, 
cf. \cite{wall} Theorem 3.5.5, 
\begin{equation}\label{baP}
  \begin{array}[t] {lccl}
(m; \beta_1 -m, \beta_2 -m, \ldots, \beta_q-m ) & & 
\text { if } & \beta_1 > 2m \\
(\beta_1 - m; m, \beta_2 -\beta_1 + m, \ldots, \beta_q-\beta_1 +m ) & &
\text { if } & \beta_1 < 2m, \, (\beta_1 -m) \! \! \not | \, m \\
(\beta_1 - m; \beta_2 -\beta_1 +m, \ldots, \beta_q-\beta_1 +m ) & & \text { if } & 
\beta_1 < 2m, \, (\beta_1 -m) | \, m 
\end{array} 
\end{equation} 
Consider $(\R^2,0)$ oriented.  Put an orientation on $(\wR , \tilde 0)$ 
so that $\pi$ preserves orientation at the points of $\tilde \gamma$.  
(This orientation depends on the tangent direction of  
$\gamma$ at $0$,  if we take consider the demi-branches with the
opposite tangent directions we get the opposite orientations  
on $(\wR , \tilde 0)$). 
% If we follow the computation of the proof of 
%cf. \cite{wall} Theorem 3.5.5 we get the following.  

\medskip
\begin{prop}\label{bacoeff}
The signs of the characteristic coefficients $\tilde A_i$ of $\tilde \gamma$ 
are :  
\begin{enumerate}
%\item [] 
\item[{ }(i)] 
\emph{Case 1 of \eqref{baP}: } $\sgn \tilde A_i = \sgn A_i$, 
$\quad i=1,\ldots, q$.  
\item[{ }(ii)]  
\emph{Case 2 of \eqref{baP}: } $\sgn \tilde A_1 = - \sgn A_1$ and 
$\sgn \tilde A_i = \sgn A_i$,  $ i=2,\ldots, q$.  
\item[{ } (iii)] 
\emph{Case 3 of \eqref{baP}: } 
$\sgn \tilde A_i = \sgn A_{i+1}$, $\quad i=1,\ldots, q-1$.  
\end{enumerate}
\end{prop} 

\begin{proof}
We may choose a system of coordinates $x,y$ at $0\in \R^2$ so that  
$\gamma(t) = (\lambda (t^m), t^m)$ and $\lambda$ is given by \eqref{fps} with 
$m_1 = \beta_1$.  Then $\tilde x := x/y$, $\tilde y :=y$, is 
a system of coordinates at $\tilde 0$ on $\wR$ and $\tilde \gamma$ is given by 
$\tilde \gamma (t) = (\tilde \lambda (t^m), t^m) $, where 
\begin{eqnarray*}
\tilde \lambda (\tilde y) = a_{m_1 / m}  \tilde y^ {(\beta_1-m) / m} + a_{m_2 / m}  
\tilde y^{(m_2-m) / m}
+ \cdots ,  \qquad \tilde y\ge 0
\end{eqnarray*}

In Case 1 of \eqref{baP}, $ {(\beta_1-m) / m}>1$ and the characteristic coefficients 
satisfy $\tilde A_i = A_i$.  

In Case 2 of \eqref{baP} we reparametrise $\tilde \gamma$  as follows.  Let  
$$
\tilde \lambda (t^m) = A_1 t^{\beta_1 - m} + \cdots = A_1 (\tilde
t(t))^{\beta_1 - m}, \quad 
  \tilde t(t) = t + \cdots .   
$$
 In the coordinates $\tilde y, \tilde x/A_1$ 
$$
\tilde \gamma (\tilde t) = ( \delta (\tilde t^{\beta_1 - m}), \tilde t^{\beta_1 - m}),
$$ 
where in $\delta (\tilde t^{\beta_1 - m}) = (t(\tilde t))^m =  \tilde t^m + \cdots$.  

The following lemma can be obtained by direct computation. 

\begin{lem}
Let the series $ t(\tilde t) = \tilde t + \sum_{k=2}^\infty b_k \tilde t^k$ be defined by 
\begin{equation}
 t^{\beta_1 - m} + \sum_{j = \beta_1 +1}^\infty a_j t^{j-m} = (\tilde t(t))^{\beta_1 - m}.
\end{equation}
Then  for each $k\ge 2$
$$
b_k = \frac {- a_{k+\beta_1-1}}{(\beta_1-m)} 
+ P_k( a_{\beta_1+1}, \ldots , a_{k+\beta_1-2} ), 
$$
where $P_k$ is a polynomial.  Moreover, if the coefficient  at 
$(a_{\beta_1+j_1})^{\alpha_1} \cdots  (a_{\beta_1+j_r})^{\alpha_r}$ in
$P_k$ is 
non-zero then $k-1 = \alpha_1 j_1 + \cdots +  \alpha_r j_r $. \qed
\end{lem}
\smallskip

By lemma the coefficient at $\tilde t^{\beta_i-\beta_1+1}$, 
$i>1$, in $ t(\tilde t) $  equals 
$\frac { - A_i }{(\beta_1-m) A_1}$.  
Consequently the coefficient at 
$\tilde t^{\beta_i-\beta_1 +m}$, $i>1$, in $\delta (\tilde t^{\beta_1 - m})$ equals 
$$
\tilde A_i = \frac { - m A_i }{(\beta_1-m) A_1}.
$$

The coordinates $\tilde y, \tilde x/A_1$ give the chosen orientation on 
$(\wR , \tilde 0)$ iff $A_1<0$.  This shows Proposition \ref{bacoeff} Case 2.  

The proof of Case 3 of \eqref{baP} is similar.  
\end{proof}

\begin{prop}\label{orders}
Let $\gamma_1, \gamma_2$ be reduced real analytic demi-branches 
tangent at $0$ and denote by $\tilde \gamma_1$ and $\tilde \gamma_2$ 
their  strict transforms by the blowing-up of the origin. 

Denote by $m$ the multiplicity of $\gamma_1$, and if $m>1$, i.e. if
$\gamma_1$ is not smooth, by $\beta_1$ the second 
exponent of the characteristic sequence of $\gamma_1$.  
The corresponding exponents of $\gamma_2$ 
will be denoted by $n$ and $\alpha_1$. Then
\begin{enumerate}
\item
If $\beta_1/m >2$ or $m=1$, and  $\alpha_1/n >2$ or $n=1$, then 
$$
O(\tilde \gamma_1,\tilde \gamma_2) = O(\gamma_1,\gamma_2) -1.
$$
\item
If $\beta_1/m < 2$, and  $\alpha_1/n >2$ or $n=1$, then 
$$
O(\tilde \gamma_1,\tilde \gamma_2) = 1.
$$
\item
If $\beta_1/m<2$, $\alpha_1/n<2$ and the first characteristic 
coefficients of $\gamma_1$ and $\gamma_2$ are of opposite sign then 
$$
O(\tilde \gamma_1,\tilde \gamma_2) = 1.
$$
\item
If $\beta_1/m < \alpha_1/n <2$ and the first characteristic 
coefficients of $\gamma_1$ and $\gamma_2$ are of same sign then 
$$
O(\tilde \gamma_1,\tilde \gamma_2) 
= \frac n {\alpha_1-n} 
$$ 
\item
If $\beta_1/m = \alpha_1/n <2$ and the first characteristic 
coefficients of $\gamma_1$ and $\gamma_2$ are of same sign then 
$$
O(\tilde \gamma_1,\tilde \gamma_2) =  
\frac m {\beta_1 -m} \, O(\gamma_1,\gamma_2) -1.
$$
\end{enumerate}
\end{prop}

\begin{proof}[Proof of Proposition \ref{orders}] 
Choose a system of coordinates so that 
$\gamma_1 (t) = (\lambda_1 (t^m), t^m)$,  
$\gamma_2 (t) = (\lambda_2 (t^n), t^n)$,  
where $\lambda_1$ and $\lambda_2$ are as in 
\eqref{fps}. Then Cases (1) and (2) are easy.  In Case (1) all exponents in 
$\lambda_1$ and $\lambda_2$ are at least $2$ and 
the claim easily follows.  
In Case (2) $\tilde \gamma_1$ and $\tilde \gamma_2$ are not tangent.   

Suppose that $\beta_1/m<2$ and $\alpha_1/n<2$ and write 
\begin{eqnarray*}
& & \lambda_1 (y) = A_1 y^{\beta_1/m} + \cdots ,\\
& & \lambda_2 (y) = B_1 y^{\alpha_1/n} + \cdots .
\end{eqnarray*}
Then both $\tilde \gamma_1$ and $\tilde \gamma_2$ are tangent to the exceptional divisor.  
If the signs of $A_1$ and $B_1$ are opposite then $\tilde \gamma_1$ and 
$\tilde \gamma_2$ have opposite tangent directions.  This shows (3).  

Suppose now that $A_1>0$ and $B_1>0$.   
In order to compare $\tilde \gamma_1$ and $\tilde \gamma_2$, as in the proof of 
Proposition \ref{bacoeff},  we have to 
pass to the coordinates $\tilde y :=y, \tilde x := x/y$, and reparametrise the 
demi-branches so we  express $\tilde \gamma_i$ in terms of $\delta_i(\tilde x)$, $i=1,2$. 
The leading exponent of $\delta_1(\tilde x)$, resp. $\delta_2(\tilde x)$, is 
$\frac m {\beta_1 - m}$, resp. $\frac n {\alpha_1 - n}$.  This shows (4) if   
$\beta_1/m<\alpha_1/n$.  

It remains to consider the case $\beta_1/m=\alpha_1/n$.  If $A_1\ne B_1$, i.e. 
$\beta_1/m = \alpha_1/n = O(\gamma_1,\gamma_2)$, then 
the leading coefficients of $\delta_1(\tilde x)$ and $\delta_2(\tilde x)$ are 
different.  This completes the remaining case of (4).  If $A_1= B_1$ then, 
as in the proof of 
Proposition \ref{bacoeff},  we use the coordinates $\tilde y , \tilde x/A_1$.  
Let $\xi =  O(\gamma_1,\gamma_2)$.  
Then the computation of the proof of 
Proposition \ref{bacoeff} for $\xi$ in place of $ \beta_i/m$  and shows that 
the first different coefficients are the ones at 
 $\tilde x^{(\beta_i - \beta_1 +m)/(\beta_1 -m)}$ as 
claimed in (5). 
\end{proof}

%%%%%%%%%%%%%%%%%%%%%%%%%%%%%%%%%%%%%%%%%%%%%%%%%%%%%%%%%%%%%%%
%\smallskip
\subsection{The symmetric demi-branch.}\label{symmetric}

%A real analytic arc has two demi-branches.  
Let $\gamma(t) = (\lambda (t^m), t^m), t\ge0$, be a reduced real analytic 
demi-branch, where  $\lambda (y)$ is given by \eqref{Ppairs}. 
\emph{The symmetric demi-branch $\gamma_-$ of $\gamma$} is obtained 
by replacing $t$ by $-t$.  If $m$ is odd it corresponds to replacing $y$ 
by $-y$ in \eqref{fps}.  Both demi-branches have the same Puiseux characteristic 
sequences and if $m$ is odd the signs of  characteristic coefficients are also 
the same.  

Suppose now that $m$ is even.  Then $\gamma_-$ and $\gamma$ are
tangent and we may compare their forms \eqref{fps} in the same system
of coordinates.  Then   
\begin{eqnarray}\label{sPpairs} 
 \qquad  \qquad  \qquad   \lambda_- (y)  =   \sum_j a_{j/m} (-1)^j y^{j/m} .
\end{eqnarray}
Let $m=2^s(2l+1)$, $s>0$, and let $d_p$ be the 
last even number in the sequence 
$d_1, \ldots, d_q$ (thus $\beta_p, \ldots, \beta_q$ are odd).  
Then the characteristic coefficients of $\gamma_-$ are 
$A_i (\gamma )$ if $i<p$ and -$A_i (\gamma )$ if $i\ge p$.  The order
contact between both demi-branches is  
$O(\gamma,\gamma_-) = \beta_p/m$.

%%%%%%%%%%%%%%%%%%%%%%%%%%%%%%%%%%%%%%%%%%%%%%%%%%%%%%%%%%%%%%%%%%%%%%%%

%%%%%%%%%%%%%%%%%%%%%%%%%%%%%%%%%%%%%%%%%%%%%%%%%%%%%%%%%%%%%%%%%%%
\medskip
\section{Cascade blow-analytic homeomorphisms}
\label{Cascade}
\smallskip

Recall that a blow-analytic homeomorphism $h: (\R^2,0) \to (\R^2,0)$ 
is called \emph{cascade} if there exists a commutative diagram 
\begin{equation}\label{cascade2}
\minCDarrowwidth 1pt
\begin{CD}
M @> b_{k}>> M_{k-1}@> b_{k-1}>> \cdots @> b_{2}>> M_1  @> b_{1} >> \R^2 \\
@V\Phi VV  @ VVh_{k-1}V @. @ VVh_{1}V @ VVhV  \\
M' @> b'_{k}>>  M'_{k-1}@> b'_{k-1}>> 
\cdots @> b'_{2}>>  M'_1 
 @> b'_{1} >> \R^2  ,
\end{CD}
\end{equation}
where $b_i$ and $ b'_i$ are point blowings-up, 
$h_i$ are homeomorphisms induced by an analytic isomorphism $\Phi$.

%Since $\varphi$ is a blow-analytic homeorphism it sends 
%a real analytic arc $\gamma :(\R,0) \to (\R^2,0) $ to a 
%real analytic arc $h\circ \gamma$.  

If a real analytic arc  $\gamma :(\R,0) \to (\R^2,0) $ 
is injective then its image is an 
 irreducible  real analytic curve 
 germ $(X,0) \subset  (\R^2,0)$ of dimension $1$.  
Let $(X_1,0)$, $(X_2,0)$ be two such curve   
germs.  If $(X_1,0)$, $(X_2,0)$ are distinct then we define 
\emph{the intersection number} 
$(X_1, X_2)_0 \in \N$ by taking the intersection number 
of the complexifications.  
%If $(X_1,0)=(X_2,0)$ then we  define $(X_1, X_1)_0 = 0$ by convention.  
By \emph{a real analytic 1-cycle at $0\in \R^2$} we mean a formal 
sum with integer coefficients of such curve germs.  
Thus, by additivity we may consider 
the intersection number at $0$ between two real analytic 1-cycles 
whose supports do intersect only at the origin. 

%If $\gamma (t)$ is reduced real analytic demi-branch and $h:(\R^2,0) \to (\R^2,0)$ is a blow-analytic
% (not necessarily cascade) homeomorphism then $h\circ \gamma$ is again reduced.  

\begin{thm}\label{cascadeproperties}
Let $h: (\R^2,0) \to (\R^2,0)$ be a cascade blow-analytic homeomorphism. 
Then $h$ preserves:
\begin{enumerate}
%\item [] 
\item[{ }(a)]
The tangency of real analytic arcs and real analytic demi-branches.
\item[{ }(b)]
The Puiseux characteristics sequence and hence the
multiplicity  
of a reduced real analytic demi-branch.  
\item[{ }(c)]
The order of contact between two reduced real analytic demi-branches. 
\item[{ } (d)] 
The intersection number between two real analytic 1-cycles at $0$ 
(of  supports intersecting only at the origin). 
\item[{ } (e)] 
If, moreover, $h$ is orientation preserving, then it 
preserves the signs of characteristic coefficients of a reduced real 
analytic demi-branch. 
\end{enumerate}
\end{thm}

\begin{proof}
The proof is by induction on the number $k$ of point blowings-up  in
\eqref{cascade2}. 
We assume that $h$ preserves the orientation. 
% otherwise we reverse the orientation of the source or of the target. 
(a) is obvious by the existence of $h_1$. 

Let $\gamma$ be a reduced real analytic arc and let $(X_\gamma ,0) \subset  (\R^2,0)$ 
by the image of $\gamma$.  Then the multiplicity of $\gamma$ equals the 
multiplicity of $X_\gamma$ at $0$ and this equals the intersection number 
of the strict transform of $(X_\gamma ,0)$ by $b_1$ with the exceptional divisor.  
Thus the invariance of the multiplicity follows from the inductive assumption.  

Now we show (d).  Let $X_1, X_2$ be two distinct irreducible curves, images of 
$\gamma_1, \gamma_2$.  We suppose that $X_1, X_2$ are tangent at the origin and denote 
by $ \tilde X_1, \tilde X_2$ their respective strict transforms intersecting at $\tilde 0$.  
Then  
%see e.g. \cite{wall} Lemma 4.4.1, 
$(X_1,X_2)_0 = m_0(X_1)m_0(X_2) + (\tilde X_1,\tilde X_2)_{\tilde 0}$ and (d) follows 
by induction.  

Let $\gamma$ be a reduced real analytic demi-branch and denote 
by $\tilde \gamma$ its strict transform by $\pi_1$.  Then, by \eqref{baP}, 
the Puiseux characteristics of $\gamma$ can be expressed in terms of  
the multiplicity of $\gamma$ and the Puiseux characteristics of $\tilde \gamma$.  
Thus (b) follows by the inductive assumption and the invariance of the multiplicity. 

Let $\gamma$ is a reduced real 
analytic demi-branch of multiplicity $m>1$ such that the second 
exponent $\beta_1$ of the characteristic sequence of $\gamma$ satisfies $\beta_1 <2m$.   
Let $\gamma_1$ be any smooth real analytic demi-branch tangent to $\gamma$.  
Then the first characteristic coefficients of $\gamma$ is positive 
if and only if $\gamma$ follows $\gamma_1$ in the clock-wise direction and therefore 
its sign is preserved by $h$ (thus we have shown a special case of (e)).   

Consequently (c) follows from Proposition \ref{orders}, (a), (b), 
the just proven special case 
of (e),  and the inductive assumption.  

The general case of (e) now follows from the same argument that we used to show 
Proposition \ref{cs}. 
\end{proof}

\begin{rem}  
Let $f_t : (\R^3,0) \to  (\R,0)$, $t \in \R$, be the Brianc\c{o}n-Speder 
family
defined by $f_t(x,y,z) = z^5 + tzy^6 + y^7x + x^{15}$.
Although $f_0$ and $f_{-1}$ are blow-analytically equivalent,
any blow-analytic homeomorphism that gives the blow-analytic
equivalence between them does not preserve the order of contact
between some analytic arcs contained in $f_0^{-1}(0)$,
see \cite{koike}. 
\end{rem}

\begin{cor}\label{cascadedistance}
Let $h: (\R^2,0) \to (\R^2,0)$ be a cascade blow-analytic homeomorphism. 
Then there exists a constant $C>0$ such that for all $(x,y)$ close to the 
origin
$$
C\inv \|(x,y)\| \le \|h(x,y)\| \le C \|(x,y)\| .
$$
\end{cor}

\begin{proof}
This follows from the invariance of multiplicity.  Indeed, 
by the curve selection lemma, it suffices to check it on a real 
analytic demi-branch $\gamma (t)$ and we 
 may assume that $\gamma$ is reduced.  Then 
$\| \gamma (t)\| $ is of size $t^m$, where $m$ denotes the 
multiplicity of $\gamma$. 
\end{proof}

We also note the following obvious property of cascade blow-analytic 
homeomorphisms. 

\begin{prop}\label{cascadedeterminant}
Let $h: (\R^2,0) \to (\R^2,0)$ be 
a cascade blow-analytic homeomorphism given by 
\eqref{cascade2}.  
Denote $\mu = b_1 \circ b_2 \circ \cdots \circ b_k$ and 
$\mu' = b'_1 \circ b'_2 \circ \cdots \circ b'_k$.  
Then the Jacobian determinants 
of $\mu'\circ \Phi$ and $\mu$, 
defined in any local system of coordinates on $M$, 
are equal up to a multiplication by a unit. 
As a consequence there exists a constant 
$C>0$ such that for  $(x,y)\ne (0,0)$ and close to the 
origin
$$
C\inv \le Jac(h) (x,y) \le C.
$$\qed
\end{prop}

%%%%%%%%%%%%%%%%%%%%%%%%%%%%%%%%%%%%%%%%%%%%%%%%%%%%%%%%%%%%%%%%%%%%%%%

%%%%%%%%%%%%%%%%%%%%%%%%%%%%%%%%%%%%%%%%%%%%%%%%%%%%%%%%%%%%%%%%%%

\section{Real tree model}
\label{realtreemodel}

%%%%%%%%%%%%%%%%%%%%%%%%%%%%%%%%%%%%%%%%%%%%%%%%%%%%%%%%%%%%%%%

\subsection{Tree model of a two variable complex analytic function germ}

Let $f(x,y)$ be a complex analytic function germ of multplicity $m$ and 
\emph{mini-regular in $x$}, that is 
$$
f(x,y) = u(x,y) (x^m + \sum_{i=1}^m a_i(y) x^{m-i}),
$$
where $m=\mult_0 f$, $u, a_i$ are analytic and $u(0,0)\ne 0$.  
Let $ x=\lambda_i (y), \, i=1,\ldots,m$, be the complex Newton-Puiseux roots of $f$.   Define \emph{the contact order of}  
$\lambda_i$ and $\lambda_j$ as  
$$
O(\lambda_i,\lambda_j) := \ord_0 \, (\lambda_i - \lambda_j)(y).
$$
Let $h\in \Q$.  We say that $\lambda_i, \lambda_j$ are \emph{congruent modulo}
 $h^+$ if $O(\lambda_i,\lambda_j)>h$.  

The tree model $T(f)$ of $f$ is defined as follows, see \cite{kuolu}
for details. 
First, draw a vertical line segment as the {\it main trunk} 
of the tree.  Mark $m=\mult_0 f(x,y)$ alongside the trunk to indicate that 
$m$ roots are bundled together.

Let $h_0:= \min \{O(\lambda_i,\lambda_j) | 1\le i,j\le m\}$. Then draw a bar, 
$B_0$, on top of the main trunk. Call $h(B_0):= h_0$ the
{\it height} of $B_0$ and mark it on the tree.  

The roots are divided into equivalence
classes modulo $h_0^+$.  We then represent each equivalence class
by a vertical line segment drawn on top of $B_0$.  Each is called
a {\it trunk}.  If a trunk consists of $s$ roots  we say it has 
{\it  multiplicity} $s$, and mark $s$ alongside (usually if $s=1$ we
do not mark it). 

Now, the same construction is repeated recursively on each trunk, getting more 
bars, then more trunks, etc..  The height of each bar and the multiplicity trunk, 
are defined likewise.  Each trunk has a unique bar on top of it.  
The construction terminates at the stage where the bar has infinite height, 
that is  on top of a trunk that contains a single, maybe multiple, root of $f$.

%%%%%%%%%%%%%%%%%%%%%%%%%%%%%%%%%%%%%%%%%%%%%%%%%%%%%%%%%%%%%%%%%%%%%%%%%%%%%%%%%%
\begin{example}\label{ex1.1}
The tree model of  $f(x,y) = (x+y)(x^2+y^3)(x^3-y^5)$. 
\vspace{4mm}
\epsfxsize=6cm
\epsfysize=4cm
$$\epsfbox{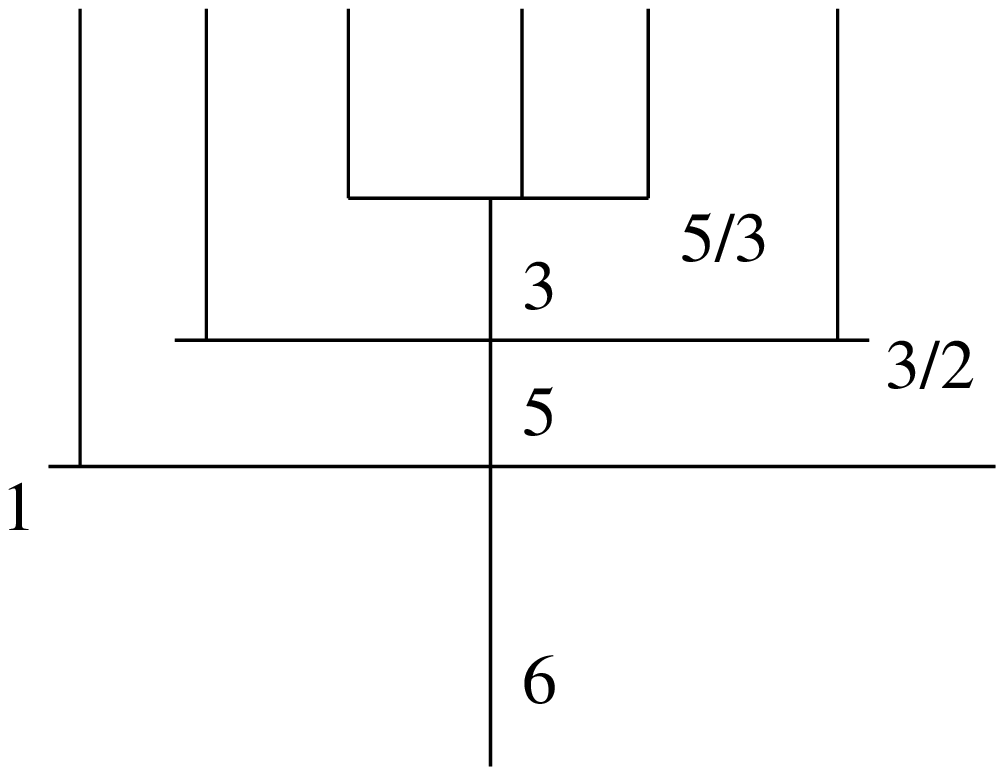}$$
\end{example}

%%%%%%%%%%%%%%%%%%%%%%%%%%%%%%%%%%%%%%%%%%%%%%%%%%%%%%%%%%%%%%%%%%%%%%%%%%%%%%%%%%

The sets of roots corresponding to trunks are called in 
\cite {kurdykapaunescu} 
\emph{bunches}. Thus each bunch $A$ is the set of roots going through
a unique bar 
$B(A)$,  one may say in this case that $A$ is the bunch bounded by $B(A)$. 
In this way we establish a one-to-one correspondence between trunks, bars, 
and bunches of roots.    

Fix a bunch $A$, $B=B(A)$, with finite height denoted $h(A)$ or $h(B)$.  
Take a root $\lambda_i (y)$ of $A$.  Let $\lambda_A(y)$ denote $\lambda _i(y)$ 
with all terms $y^e$, $e\ge h$, omitted. 
%(In particular, $\lambda_{B_*}(y)=0$.)
Clearly, $\lambda_A$ depends only on $A$, not on the choice of $\lambda_i\in A$.  
We can then write for each $\lambda_i (y) \in A$ 
\begin{equation}\label{lambdaA}
\lambda_{i}(y) = \lambda_A (y) + c_i y^{h(A)} +  \cdots , \qquad c_i\in \C.
\end{equation}

\smallskip
\begin{rem}
%Fix a Newton-Puiseux root $\lambda$.  Then the denominators in the expansion of% $\lambda$ 
%divide the commom denominator of the heights of all the bars that bo%und $\lambda$.
The fractional power series $\lambda_i (y)$ are  well-defined only up to 
the action of a group of roots of unity.  
One may make them well-defined by fixing the argument of $y\in \C$, 
for instance  
$y\in \R, y>0$, cf. \cite{wall} p. 98.  

Fix a bunch $A$.  
Suppose that the denominator of $h(A)$ does not divide the common denominator 
of exponents of $\lambda_A$.  It means that these roots of $A$ for which 
$c_i\ne 0$ get a new Puiseux pair at the bar $B(A)$ 
and those with $c_i= 0$ do not. Formally, for a particular root, 
 the tree does not contain the  information what is   
the coefficient $c_i$ at $y^{h(A)} $ or even  whether $c_i\ne0$.  
Nevertheless, by counting the number of subbunches of $A$, one can read from the tree 
 whether there exists a sub-bunch that does not take a new 
Puiseux pair at $B(A)$.  If it exists it must be unique.   
\end{rem}

%%%%%%%%%%%%%%%%%%%%%%%%%%%%%%%%%%%%%%%%%%%

%\smallskip
\subsection{Real part of tree model}

Suppose that $f(x,y)$ is real analytic.  
Consider the Newton-Puiseux roots as arcs $x=\lambda_i(y)$ 
defined for $y\in \R, y\ge 0$. 
%Identify (arbitrarily) the roots to the terminal 
%vertical lines on the tree.  
The complex conjugation acts on the 
roots, and hence on the tree model $T(f)$.
% by conjugating the coefficients of \eqref{lambdaA}.  
A bunch $A$ of $T (f)$ is called \emph {real} 
if it is stable by complex conjugation, or equivalently if $\lambda_A$ is real. 
A bar or a trunk is \emph{real} if and only if so is the corresponding bunch.  
After \cite {kurdykapaunescu} the conjugation invariant part of $T(f)$, 
that we denote by $T_+(f)$, is called \emph {the real part of $T(f)$}.  

Similarly, by fixing $y\le 0$, we may define $T_-(f)$.  We may identify 
$T_+(f)$ and $T_-(f)$ if all denominators of the exponents of \eqref{lambdaA} 
are odd.  But in general $T_+(f)$ and $T_-(f)$ are different.

%%%%%%%%%%%%%%%%%%%%%%%%%%%%%%%%%%%%%%%%%%%%%%%%%%%%%%%%%%%%%%%%%%%%%%%
\begin{example}%\label{ex1.1}
We draw $T_+(f)$ and $T_-(f)$ for   $f(x,y) = (x+y)(x^2+y^3)(x^3-y^5)$
of Example \ref {ex1.1}
%\vspace{3mm}
\epsfxsize=12cm
\epsfysize=4cm
$$\epsfbox{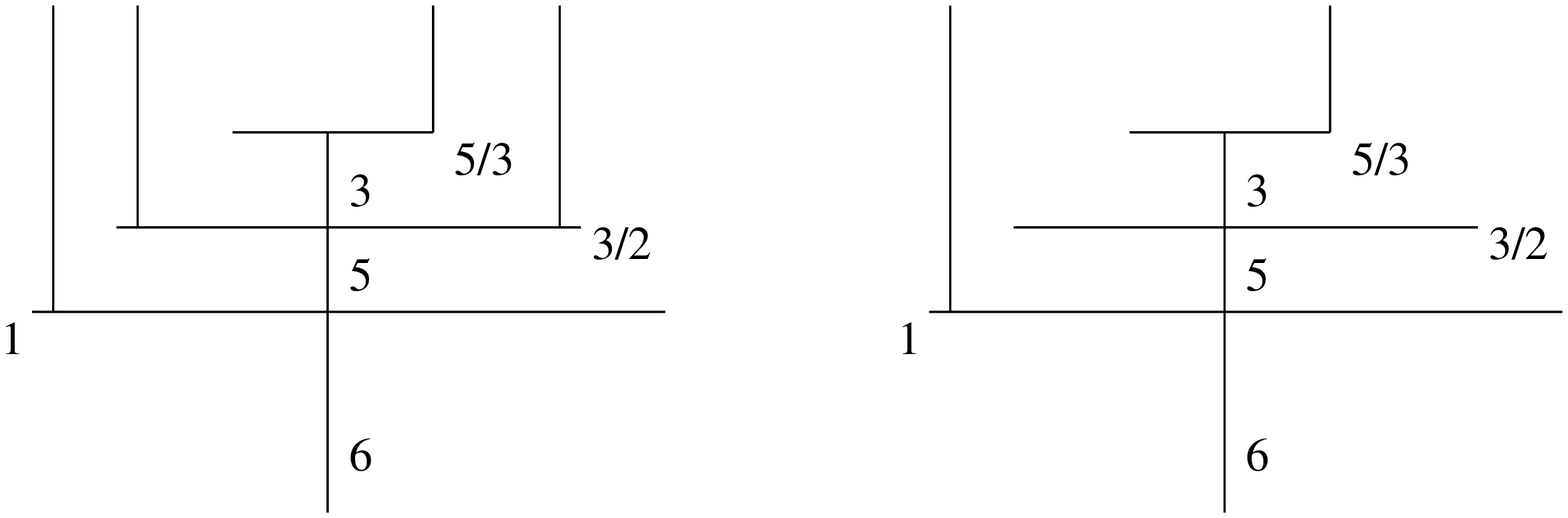}$$
\smallskip
\centerline{$T_- (f)$\hspace{55mm} $T_+ (f)$\hspace{5mm}}
\end{example}
%%%%%%%%%%%%%%%%%%%%%%%%%%%%%%%%%%%%%%%%%%%%%%%%%%%%%%%%%%%%%%%%%%%%%% \smallskip

\begin{rem}
In \cite {kurdykapaunescu} the authors define also the almost real bunches.
For $f$ real analytic the notions of a real bunch and of an almost real bunch 
coincide.
\end{rem}

\begin{rem}\label{complicated}
The real part of tree model, even if all denominators in $\lambda_i$ 
 are odd, does not determine the Newton-Puiseux pairs 
of the real roots.  Indeed,   
Let $f(x,y) = x(x^3 - y^5)((x^3 - y^5)^3 - y^{17})$, 
$g(x,y) = x(x^3 + y^5)(x^3 - y^7)(x^6 + x^3y^5 + y^{10})$.  All singular roots of $f$ 
have a single Puiseux pair $(5,3)$ and two of such roots have contact order $7/3$.  
The singular roots of $g$ have also a single Puiseux pair, either $(5,3)$ or $(7,3)$.   
Then the real parts $T_+(f)=T_+(g)$ of the tree models coincide, but
the real tree models of $f$ and $g$, see next subsection, are
different. As a consequence of Theorem \ref{allequivalent} we see that
$f$ and $g$  are not blow-analytically equivalent.  See also example \ref{complicatedex}

%%%%%%%%%%%%%%%%%%%%%%%%%%%%%%%%%%%%%%%%%%%%%%%%%%%%%%%%%%
\vspace{3mm}
\epsfxsize=5cm
\epsfysize=3.2cm
$$\epsfbox{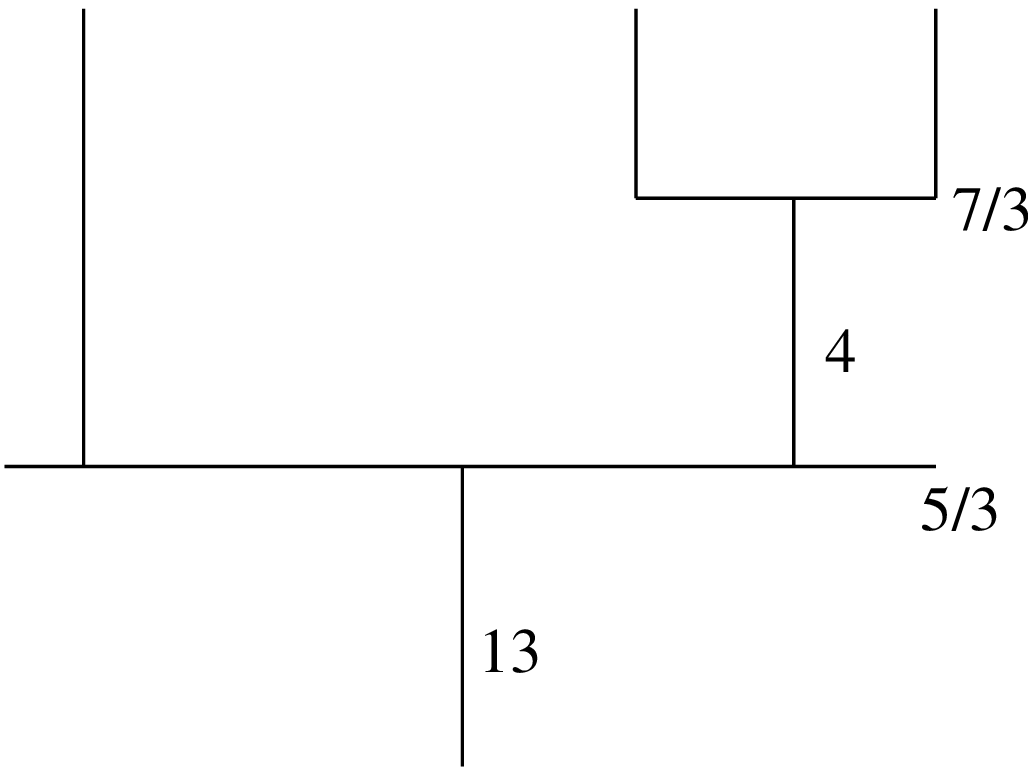}$$
%\smallskip
%\centerline{$T_+(f)$\hspace{55mm} $T_-(f)$\hspace{5mm}}
\end{rem}

%%%%%%%%%%%%%%%%%%%%%%%%%%%%%%%%%%%%%%%%%%%

\smallskip
\subsection{Real tree model of $f$}

Let $f : (\R^2,0) \to (\R,0)$ be a real analytic function germ.  Fix $v$ a unit 
vector of $\R^2$.  Fix any local system of coordinates $x,y$ such that \par
\emph{- $f(x,y)$ is mini-regular in $x$}\par
\emph{- $v$ is of the form  $(v_1,v_2)$ with $v_2 >0$.}\\
Consider the Newton-Puiseux roots as arcs $x=\lambda_i(y)$ 
defined for $y\in \R, y\ge 0$. 

\smallskip
We define \emph{the real tree model of $f$ relative to $v$}, and denote it by 
$\R T_v (f)$, as the part of $T_+(f)$ consisting only of the roots tangent to $v$ 
with the following additional information.  Let $A$ be a real bunch such that 
$B=B(A)$ is a bar of $\R T_v (f)$.  Then : \\
{- draw the trunks on $B$ realising the subbunches of $A$ keeping the clockwise 
order of the roots (i.e. the order of the coefficients $c_i$ in \eqref{lambdaA}), \\
{- whenever $B$ gives a new Puiseux pair to some roots of $A$  
(that can be easily computed from $T_+(f)$) we mark $0$ on $B$ and grow at it 
the unique sub-bunch of $A$ with $c_i=0$, i.e. consisting of the roots that 
do not have the new Puiseux pair at $B$.  Hence we are able to determine from 
the tree also the sub-bunches with positive and negative $c_i$.  
Graphically,  we identify $0\in B$ with the point of $B$ 
that belongs to the trunk  supporting $B$.

%%%%%%%%%%%%%%%%%%%%%%%%%%%%%%%%%%%%%%%%%%%%%%%%%%%%%%%%%%%%%%%%%%%%%%%
\begin{example}\label{tree}
Consider $f(x,y) = x(x^3 - y^5)(x^3 + y^5)$ and 
$g(x,y) = x(x^3 - y^5)(x^3 - 2y^5)$ of Example \ref{examplesimple}.  
The real trees $\R T_{(0,1)}(f),$ and $\R T_{(0,1)}(g)$, drawn below, 
 are different.  
\vspace{5mm}
\epsfxsize=10cm
\epsfysize=3cm
$$\epsfbox{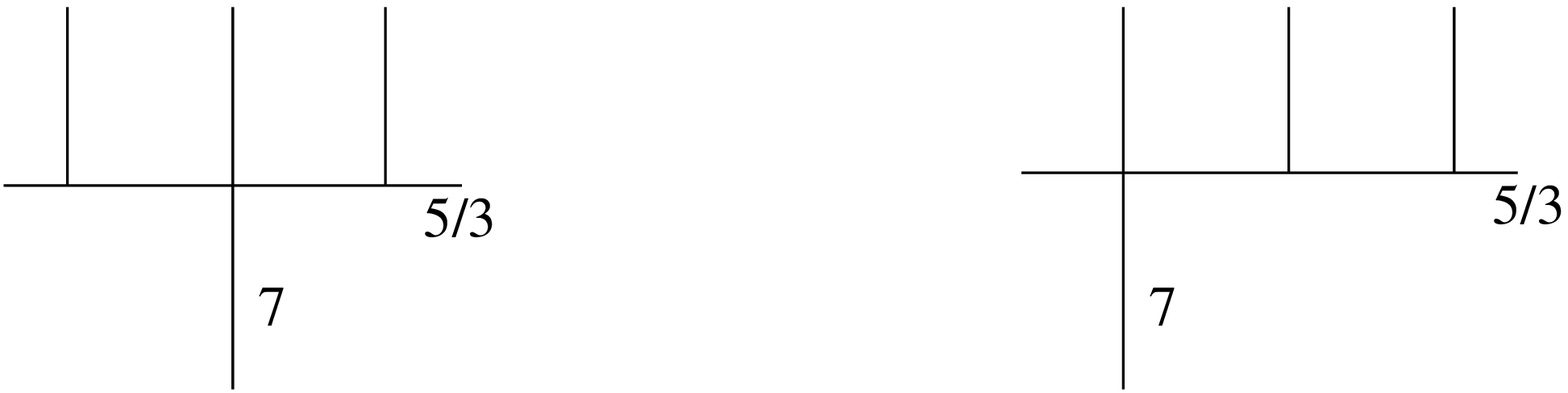}$$
\centerline{\hspace{10mm}$\R T_{(0,1)}(f)$\hspace{80mm} $\R T_{(0,1)}(g)$\hspace{10mm}}
\vspace{3mm}
\end{example}
%%%%%%%%%%%%%%%%%%%%%%%%%%%%%%%%%%%%%%%%%%%%%%%%%%%%%%%%%%%%%%%%%%%%%%%%%%

\begin{example}\label{complicatedex}
Consider  $f(x,y) = x(x^3 - y^5)((x^3 - y^5)^3 - y^{17})$, 
$g(x,y) = x(x^3 + y^5)(x^3 - y^7)(x^6 + x^3y^5 + y^{10})$ as in remark \ref{complicated}. 
The real trees $\R T_{(0,1)}(f),$ and $\R T_{(0,1)}(g)$, drawn below, 
 are different.  

%%%%%%%%%%%%%%%%%%%%%%%%%%%%%%%%%%%%%%%%%%%%%%%%%%%%%%%%%%
\vspace{3mm}
\epsfxsize=10cm
\epsfysize=3.2cm
$$\epsfbox{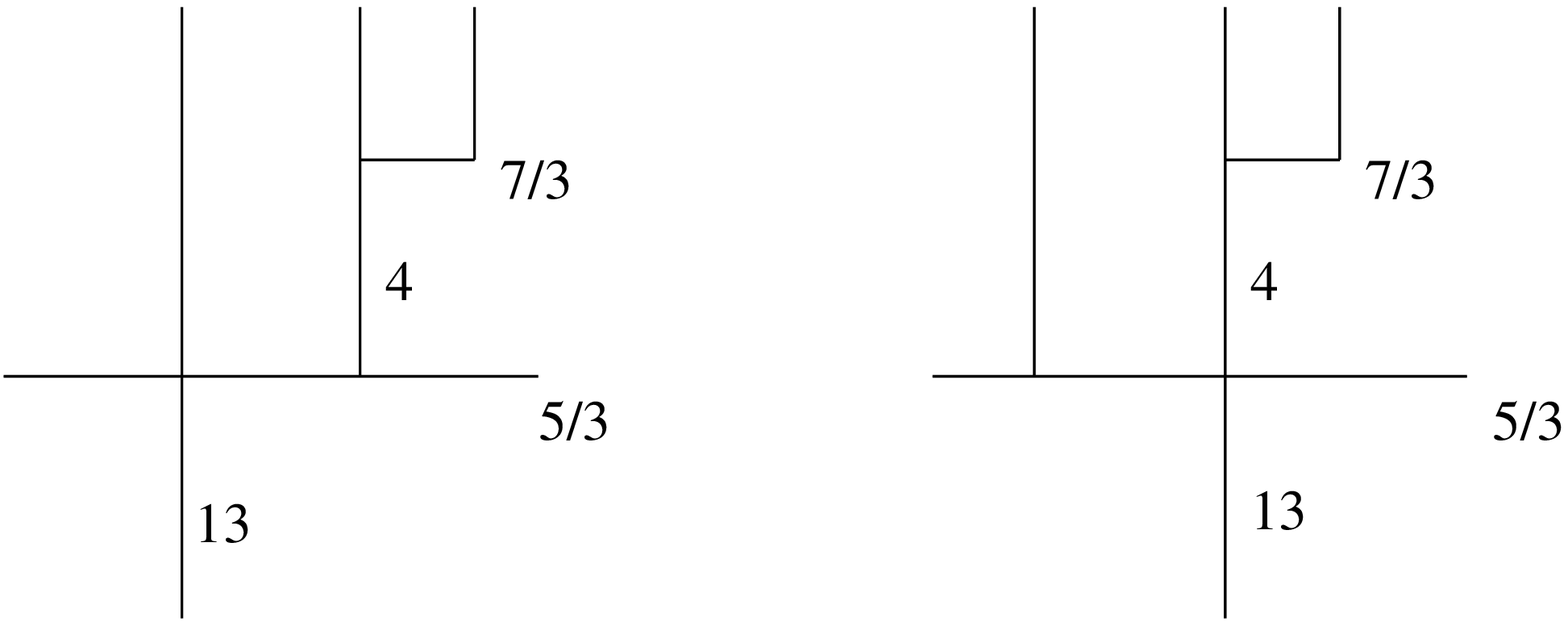}$$
%\smallskip
%\centerline{$T_+(f)$\hspace{55mm} $T_-(f)$\hspace{5mm}}
\end{example}
\medskip

\begin{defn}
Let $f(x,y)$ be an analytic function germ.  The \emph{real tree model 
$\R T(f)$ of $f$} is defined as follow. 
\begin{itemize}
%\item 
%Draw a vertical line segment on $B_*$ as the {\it main trunk} of the tree.  
%Mark $m=\mult_0 f(x,y)$ alongside.  
\item 
Draw a bar $B_0$ that is identified with $S^1$.  
 We define $h(B_0)=1$ and call $B_0$ \emph{the ground bar}.  We mark  
$m(B_0) : = 2 \mult_0 f(x,y)$ below the ground bar.  
\item
Grow on $B_0$ non-trivial individual $\R T_v(f)$ for $v\in S^1$, keeping the 
clockwise order.
\item
Let $v_1, v_2$ be any two subsequent unit vectors for which $\R T_v(f)$ is nontrivial. 
Mark on $B_0$ of  $\R T(f)$  the sign of $f$ in the sector between $v_1$ and $v_2$. 
Note that one such sign determines all the other signs between two subsequent  
unit vectors for which $\R T_v(f)$ is nontrivial (passing $v$ changes
this sign if and only if $\R T_v(f)$ contains an odd number of roots.) 
\end{itemize}
If the leading homogeneous part $f_m$ of $f$ satisfies 
$f_m\inv (0) =0$ then $B_0$ is the only bar of $\R T(f)$. 
\end{defn}

To each real bunch $A$ of the real tree model  $\R T(f)$ we associate 
\emph{a generic demi-branch 
  associated to $A$},  $\gamma_{A,gen} : x=\lambda_{A,gen} (y)$, 
where $\lambda_{A,gen}$ is given by 
\begin{equation}\label{lambdagen}
\lambda_{A, gen}(y) = \lambda_A(y) + c y^{h(A)} +  \cdots .
\end{equation} 
 with the constant $c$
generic in $\R$.  \emph{The characteristic exponents} and respectively
\emph{the
  signs of characteristic coefficients associated to $A$} are those of 
 $\gamma_{A,gen}$ that are $\le h(A)$, respectively correspond to the
 exponents $<h(A)$.

Let $B$ be a bar of $\R T(f)$ and let $T_i$ be the trunks grown on
$B$.  Then, in the complex case, the multiplicity $m(B)$ of $B$ equals the 
sum of the multiplicities of $T_i$.  For the real tree model this is 
no longer true, $m(B)- \sum_i m(T_i)$ can be strictly positive, though it is 
always even.  
Note that in the real case each complex root is counted twice, 
once for $y\ge 0$ and once for $y\le 0$.    

\begin{example}\label{weighted} 
Let $f(x,y)$ be a singular germ  weighted homogeneous with weights
$w_1> w_2$.  Then $\R T(f)$ has exactly two trunks grown on $B_0$
corresponding to the $y$- axis, $y\le 0$ and $y\ge 0$.  Each of these
trunks is bounded by a bar of height $w_1/w_2$.  These are all bars of 
the tree.    
\end{example}

\begin{rem} \label{opposite}
Let $B=B(A)$ be  a bar of $\R T_v(f)$ such that $m, \beta_1$, $\beta_1/m<2,$ is
 the Puiseux characteristic sequence of 
$\gamma=\gamma_{A,gen} : x=\lambda_{A,gen} (y)$, 
 and $h(A) = \beta_1/m$.  Then $B$ admits 
\emph{the opposite bar} $-B$ of $\R T_{-v}(f)$ of the same height 
with the generic arc that is given  in the system of coordinates 
$x'=-x, y'=-y$ by the same formula $x'=\lambda_{-A,gen} (y')
= \lambda_{A,gen} (y')$ and $m(-B) = m(B)$.  
Moreover, let $T_0$ be the trunk that grows on $B$ at $0$, that is consisting of the roots 
that do not take a new Puiseux pair at $h(B)$.  Then there is a trunk 
that grows on $-B$ at $0$ that have the same multiplicity as $T_0$. 
\end{rem}

\begin{defn}
We call two real trees $\R T(f), \R T(g)$ \emph{isomorphic} if there is a homeomorphism
$\varphi$ of their ground bars sending one tree to the other and preserving  the 
 multiplicities and heigths of bars and signs of the characteristic coefficients.  
If, moreover, $\varphi$ preserves the orientation we call the trees \emph{orientably
  isomorphic}.  
\end{defn}

\begin{example}\label{Brieskornoriented}
Let $1<p<q$ be odd numbers.  Then the real trees of 
$f(x,y)= x^p-y^q$ and $g(x,y)= x^p +y^q$ are not orientably 
isomorphic. 
\end{example}

\medskip
%%%%%%%%%%%%%%%%%%%%%%%%%%%%%%%%%%%%%%%%%%%%%%%%%%%%%%%%%%%%%%%
\subsection{Effect of a blowing-up.}\label{blowtree}

Let $\pi: \wR  \to \R^2$ be the blowing-up of the origin. Fix a
point 
$\tilde 0\in \wR$ on the exceptional divisor such that 
$f\circ \pi$ at $\tilde 0$ 
is not the $m$-th power of the equation defining the exceptional divisor.  
Let $v$ and $-v$ be the two opposite unit vectors of $\R^2$ corresponding to 
$\tilde 0$.  Then $\R T((f\circ \pi)_{\tilde 0})$ is determined by  
$\R T_v(f)$ and $\R T_{-v}(f)$ and a choice of orientation of 
$(\wR, \tilde 0)$ in the following way.  

Fix a bar $B=B(A)$ of $\R T_v(f)$ and $\R T_{-v}(f)$. Let 
$\gamma_{A,gen}$ denote its generic demi-branch,  $h=h(A)$ its height, and 
$m(A)$ the multiplicity of $A$.  By $m, \beta_1, ...$ we denote the characteristic sequence of 
$\gamma_{A,gen}$ (note that in general $m\ne m(A)$).  
$B$ gives rise to a bar $\tilde B$ or a couple of bars $\tilde B_+, \tilde B_-$ of  
$\R T((f\circ \pi)_{\tilde 0})$ that we describe by giving their generic demi-branches,  
heights, and  multiplicities.

Denote by $\tilde \gamma$ the strict transform of 
$\gamma=\gamma_{A,gen}$.  
Then the Puiseux characteristic sequence of $\tilde \gamma$
is 
given by \eqref{baP} and the sign of characteristic exponents are
given by  
Proposition \ref{bacoeff}.  The height of new bars can be computed 
using Proposition \ref{orders}.

\begin{enumerate}
\item
If $\beta_1/m >2$ or $m=1$ ($\gamma$ smooth), then $h \ge 2$.\\  
 If $h > 2$, 
then  $\tilde \gamma$ and $h(\tilde B)  = h -1$ define a 
bar $\tilde B$, and  
$m(\tilde B) = m(B)$.\\
If $h = 2$ then $B$ gives the ground bar $\tilde B_0$ of 
$\R T((f\circ \pi)_{\tilde 0})$.
\item
The exceptional divisor is a root of  
$f\circ \pi$ of multiplicity $\mult_0 f(x,y)$.  
It gives rise to two bars of $\R T((f\circ \pi)_{\tilde 0})$ of infinite height.  
\item
Let $\beta_1/m < 2$ and $h = \beta_1/m$.  Then $B$ breaks down into
two parts relative to the sign of $c$ in \eqref{lambdagen}, i.e. 
the sign of the first characteristic
coefficient of $\gamma$.  Denote these two different generic $\gamma$ by 
$\gamma_+$ and $\gamma_-$ respectively.   \\
The strict transform  $\tilde \gamma_+$ of $\gamma_+$, 
and $h(\tilde B_+)   = \frac 
m {\beta_1 - m}$ gives  a new bar $\tilde B_+$.  One of the 
demi-branches of the exceptional divisor grows on $\tilde B_+$ and 
divides it  into two parts.  We note that $\tilde \gamma_+$ grows on 
one of these part.  The other part corresponds 
to (a half of) $-B$,  see Remark \ref{opposite}.   
Similarly the strict transform $\tilde \gamma_-$ of $\gamma_-$, 
and $h(\tilde B_-)   = \frac 
m {\beta_1 - m}$ gives  a bar $\tilde B_-$.  
%The  bar $\tilde B_-$ contains the other demi-branch of the exceptional divisor. 
The multiplicities of new bars are given by 
$$
m(\tilde B_+) = m(\tilde B_-) =  m(B) - m(A') + \mult_0 f(x,y),
$$
where $A'$ is the subbunch of $A$ of the roots that 
does not take a new Puiseux pair at $h(B)$.
% (we put $m(A')=0$ if $A'=\emptyset$).  
The strict transforms of roots of  $A'$ do not grow neither on 
$\tilde B_+$ nor on $\tilde B_-$. 
\item
Let $\beta_1/m < 2$ and $h > \beta_1/m$ ($h < \beta_1/m$ cannot
happen).  \\
Then $B$ grows over a bar $B'$  with $h(B') = \beta_1/m$ and 
gives rise to a bar  $\tilde B $ that grows either on $\tilde B'_+$ or
on $\tilde B'_-$,  $h(\tilde B) =  \frac 
m {\beta_1 - m} h - 1$ as 
follows from (5) of Propostion \ref{orders}, and 
$m(\tilde B) = \frac {\beta_1-m}m m(B)$.  
\end{enumerate}

Each bar of $\R T((f\circ \pi)_{\tilde 0})$ comes from a bar of $\R
T_v(f)$ or $\R T_{-v}(f)$.  
The only possible exception could be the ground bar $\tilde
B_0$ of $\R T((f\circ \pi)_{\tilde 0})$ if there  is no bar of 
$\R T_v(f)$ or  $\R T_{-v}(f)$ with $m(B)=1$ and $h(B) = 2$.  

The effect of a blowing-up on the tree can be also expressed in terms of horn, 
see subsection \ref{realconstruction} below.

%%%%%%%%%%%%%%%%%%%%%%%%%%%%%%%%%%%%%%%%%%%%%%%%%%%%%%%%%%%%%%%%%%%%%%%

%%%%%%%%%%%%%%%%%%%%%%%%%%%%%%%%%%%%%%%%%%%%%%%%%%%%%%%%%%%%%%%%%%%
\medskip
\section{Horns and root horns}
\label{Horns}
\smallskip

In this section we characterise the real tree model of 
$f:(\R^2,0)\to (\R,0)$ in terms of the real analytic geometry without explicit 
referring to the complex Newton-Puiseux roots of $f$.  This characterisation  
is used to show the blow-analytic invariance of the real tree model.  
The main idea is to replace a real bunch of roots by a geometric object, a horn.

%%%%%%%%%%%%%%%%%%%%%%%%%%%%%%%%%%%%%%%%%%%

\subsection{Horns}

Let $\gamma (t)$ be a reduced real analytic demi-branch given by \eqref{fps}.
Define the {\em horn-neighbourhood of $\gamma$ 
of exponent} $\xi > 1$ {\em and width} $C > 0$ by 
$$
H_{\xi}(\gamma ; C) = \{ (x,y); \, \dist ((x,y), \text {image } (\gamma)) \le
C|(x,y)|^{\xi} \} .
$$
%where $\im (\gamma)$ denotes the image of the demi-branch $\gamma$.  
By a {\em horn-neighbourhood of $\gamma$ of exponent} $\xi > 1$ we mean 
$H_{\xi}(\gamma ; C)$ for $C$ large and we denote it by $H_{\xi}(\gamma)$.  

\begin{rem}\label{C} 
In order to simplify the exposition we use the following convention 
\begin{enumerate}
\item
If $O(\gamma_1,\gamma_2)\ge \xi$ then we say that $H_{\xi}(\gamma_1)= H_{\xi}(\gamma_2)$  
by meaning that for any $C_1>0$ there is $C_2>0$  such that 
$$
H_{\xi}(\gamma_1 ; C_1) \subset H_{\xi}(\gamma_2 ; C_2), \qquad 
H_{\xi}(\gamma_2 ; C_1) \subset H_{\xi}(\gamma_1 ; C_2) 
$$
\item
For similar reasons, if the demi-branch is represented by 
 \eqref{fps}, we write   
$$
H_{\xi}(\gamma)= \{ (x,y); \ |x - \lambda (y)| \le
C|y|^{\xi} \} ,
$$
for $C$ sufficiently large.  
\end{enumerate}
\end{rem}  

Let $H = H_{\xi}(\gamma)$ be a horn.  Then by \emph {a generic arc} $\gamma_H$ 
in $H$ we mean a demi-branch given by  
\begin{equation}%\label{lambdaA}
x=\lambda_{H,gen} (y) = \lambda_(y) + c y^{\xi} +  \cdots , \qquad y\ge 0, 
\end{equation} 
where $c\in \R$ is a generic constant. The \emph {characteristic exponents of $H$} 
are those of $\gamma_H$ that are $\le \xi$.  Similarly we define \emph{the signs of 
characteristic coefficients of $H$} taking into account the exponents $<\xi$.  

Given $f$ and $\gamma$ as above.   Fix  $\xi \ge 1$ and expand 
\begin{equation}\label{genericarc}
f(\lambda (y)+zy^{\xi},y) 
= P_{f,\gamma,\xi} (z) y^{\ord_{\gamma}f(\xi )} + \cdots  , 
\end{equation}
where the dots denote higher order terms in $y$ and  $\ord _{\gamma}f(\xi)$ 
is the smallest exponent with non-zero coefficient.  
This coefficient,  
$P_{f,\gamma,\xi} (z)$,  is a  polynomial function of $z$.

%%%%%%%%%%%%%%%%%%%%%%%%%%%%%%%%%%%%%%%%%%%

\smallskip
\subsection{Root horns} 
Let $A$ be a real bunch  of $\R T_v (f)$.  
Then $A$ defines a horn 
$$
H_A := \{ (x,y); \ |x - \lambda_A (y)| \le
C|y|^{h(A)} \} ,
$$
where $C$ is a large constant.  Then $\gamma_H = \gamma_{A,gen}$, see \eqref{lambdagen}.  
A horn that equals $H_A$ for a bunch $A$ is called \emph{a root horn}.

\begin{prop}\label{barreal}
%{\rm [compare  \cite{koikeparusinski2}, Proposition 7.5] } 
Let $H $ be a horn of exponent $\xi$.  Then $H$ is a root horn  for $f(x,y)$ if and only if 
$P_{f,\gamma_H,\xi} (z) $ has at least two distinct complex roots. 
If this is the case, $H=H_A$, then $h(A) = \xi $ and $m_A = \deg P_{f,\gamma_H,\xi_H}$.
\end{prop}

\begin{proof}
Suppose that $H=H_A$ and let $A= \{\gamma_1, \ldots, \gamma_{m_B}
\} $ be the corresponding bunch of roots.  
These roots are truncations of complex Newton-Puiseux roots of $f$: 
\begin{equation}\label{complexroot}
\gamma_{\C,k}: x= \lambda _k (y) = \lambda_H(y) + a_{\xi,k} y^\xi + \cdots, \quad 1\le k\le m_A 
\end{equation}
with $\lambda_H$ real and $a_{\xi,k}\in \C$.  Denote by $\gamma_{\C,j}: x= \lambda _j (y)$, 
$j=m_A +1, \ldots, m,$ the remaining complex Newton-Puiseux roots of $f$.  Then   
\begin{equation*}%\label{genericarc}
f(\lambda_H (y)+zy^{\xi},y)  = u(x,y) \, \prod_{i=1}^m (\lambda_H (y)- \lambda_i (y)+zy^{\xi})  
= P_{f,\gamma_H,\xi} (z) y^{\ord_{\gamma_H}f(\xi )} + \cdots  , 
\end{equation*}
where $u(0,0)\ne 0$.   Note that $O(\lambda_H, \lambda_j) <\xi$ for $j>m_A$.   Therefore 
\begin{eqnarray*}
&  P_{f,\gamma_H,\xi_H} (z) & = u(0,0)  \, \prod_{i=1}^{m(A)} (z -a_{\xi ,i}),  \\
&  \ord_{\gamma_H}f(\xi_H ) & = m_A \xi_H + \sum_{j=m(A)+1}^m O(\lambda_H, \lambda_j). 
\end{eqnarray*}  
By construction of the tree there are at least two roots $\gamma_i$ and $\gamma_j$ of 
\eqref{complexroot} such that 
$O(\lambda_i, \lambda_j)  = \xi$. Thus $P_{f,\gamma_H,\xi_H} (z) $ has
at least two distinct complex roots.  

Let $H = H_{\xi}(\gamma)$ be a horn,  where 
$$
\gamma : x= \lambda_H(y) + a_\xi y^{\xi} +  \cdots. 
$$  
%\begin{equation*}%\label{genericarc}
%f(\lambda (y)+zy^{\xi},y)  = P_{f,\gamma,\xi} (z) y^{\ord_{\gamma}f(\xi )} + \cdots  . 
%\end{equation*}
By the Newton algorithm  for computing the complex Newton-Puiseux roots of $f$ 
to each root $z_0$ of $ P_{f,\gamma,\xi}$ of multiplicity $s$ correspond exactly $s$ 
Newton-Puiseux roots of $f$, counted with multiplicities, of the form  
\begin{equation*}%\label{genericlambda}
\gamma_0 : x= \lambda_H(y) + (a_\xi +z_0)  y^{\xi} +  \cdots .
\end{equation*} 
(This is essentially the way the Newton-Puiseux theorem is proved as in \cite{walker}.)  
Thus, if $ P_{f,\gamma,\xi}$ has at least two distinct roots, then there exist at least two such Newton-Puiseux 
roots with contact order equal to $\xi$.  This shows that $H$ is of the form $H_B$, as claimed.  
 \end{proof}

%%%%%%%%%%%%%%%%%%%%%%%%%%%%%%%%%%%%%%%%%%%

%\smallskip
\subsection{Real construction of the real tree model}\label{realconstruction}

The root horns can be ordered by inclusion and by the clockwise order  
around the origin.  Thus $V=H_A$ is contained in $V'=H_{A'}$ if and only if 
bunch $A$ is contained in $A'$.   
The height $h(A)$ and the multiplicity $m(A)$ are given by Proposition \ref{barreal}

Let $\gamma_{V,gen}$ be a generic arc associated to the horn $V=H_A$.  
Let $\gamma$ be any complex root in $A$.  
Then the Puiseux characteristic exponents of $\gamma$ that are 
$< h(A)$ and the corresponding   signs of characteristic coefficients are those 
of $\gamma_{V,gen}$.  In particular, let $\tilde A$ be a sub-bunch of
$A$.    
Then $\gamma_{H_{\tilde A},gen}$ contains the information whether $\tilde A$ 
takes a new Puiseux pair at $h(A)$ and, if this is the case, 
 the sign of the characteristic coefficient 
at $h(A)$.  

Thus the entire real tree model $\R T(f)$ can be obtained from the knowledge 
of the root-like horns and their numerical invariants.

\begin{cor}\label{b-agivessametrees}
Let $f:(\R^2,0)\to (\R,0)$ and $g:(\R^2,0)\to (\R,0)$ be real analytic 
function germs.  If $f$ and $g$ are blow-analytically equivalent by
an orientation preserving homeomorphism $h$ then the  
 real trees of $f$ and $g$ are orientably isomorphic.  
\end{cor}

\begin{proof}
By section \ref{lift}, $f$ and $g$ are cascade 
blow-analytically equivalent. By  Theorem \ref{cascadeproperties} and 
Corollary \ref{cascadedistance}, the image by $h$ of a root-like horn is 
a root-like horn with the same numerical invariants.  
\end{proof}

%%%%%%%%%%%%%%%%%%%%%%%%%%%%%%%%%%%%%%%%%%%%%%%%%%%%%%%%%%%%%%%%%%%%%%%%%%%%%%%%%

%%%%%%%%%%%%%%%%%%%%%%%%%%%%%%%%%%%%%%%%%%%%%%%%%%%%%%%%%%%%%%%%%%%%%%%

%%%%%%%%%%%%%%%%%%%%%%%%%%%%%%%%%%%%%%%%%%%%%%%%%%%%%%%%%%%%%%%%%%%

\section{End of Proof of Theorem \ref{allequivalent} } \label{EndofProof}

We prove the theorem only in the oriented case. In (1) of Theorem
\ref{allequivalent} it means that the blow-analytic homeomorphism preserves the
orientation.  In (2) it means that the isomorphism $\Phi : M\to \tilde
M$, cf. subsection \ref{constructing}, blows down to an orientation
preserving homeomorphism $h : (\R^2,0) \to (\R^2,0)$.  In (3) it means that  then the  
 real trees of $f$ and $g$ are orientably isomorphic.

%First we consider the oriented case.  
(1) $\Longleftrightarrow $ (2) by Proposition \ref{fukuiprop} and Section \ref{lift}. 
(2) $\Longrightarrow $ (3) by Corollary \ref{b-agivessametrees}.  
We show (3) $\Longrightarrow $ (2).  

Let $\pi: \wR  \to \R^2$ be the blowing-up of the origin. If 
both $f\circ \pi$ and $g\circ \pi$ are normal crossings then    $\pi: \wR  \to \R^2$ 
is the 
minimal resolution of $f$ and of $g$ and the claim follows. 
Otherwise, by Lemma \ref{normalcrossing}, 
there is a one-to-one correspondence between the non-normal crossings points 
of $f\circ \pi$ and $g\circ \pi$. Fix such a pair of points 
$p_1 \in \pi \inv (0)$, 
$p_2 \in \pi \inv (0)$.  By subsection \ref{blowtree} the 
real trees of $(f\circ \pi,p_1)$ and $(g\circ \pi,p_2)$, for a choice of 
local orientations, coincide.  Thus, by the inductive assumption, 
there are neighbourhoods $U_1$ and $U_2$ of $p_1$ and $p_2$ respectively, 
such that  the minimal 
resolutions of $f\circ \pi|_{U_1}$ and $g\circ \pi|_{U_2}$ are isomorphic.  
This isomorphism induces an analytic isomorphism of punctured 
neighbourhoods $U_1\setminus p_1$ and $U_2 \setminus p_2$.  
%It remains to show that  we may choose in a coherent way the local orientations at all such  pairs. 

For a pair of corresponding points  $p_1$ and $p_2$, fix a half-line 
$l_1$ in $\R^2$ at $0$ representing $p_1$. Let $l_2$ be the  half-line corresponing to 
$l_1$ by the  identification of both trees. A choice of half-line $l_1$ 
induces 
on $(\wR, p_1)$ an orientation, the one that lifts the 
canonical orientation of $\R^2$ at the points 
of $l_1$.  Note that replacing $l_1$ and $l_2$ by the opposite 
half-lines reverses both orientations, and therefore gives 
the same isomorphism of minimal resolutions of  
$(f\circ \pi,p_1)$ and $(g\circ \pi,p_2)$.  Thus  
the  identification of both trees gives a coherent system of 
orientations and  the local isomorphisms of minimal
resolutions glue together, see Remark \ref{Nash}.      \qed

%%%%%%%%%%%%%%%%%%%%%%%%%%%%%%%%%%%%%%%%%%%%%%%%%%%%%%%%%%%%%%%%%%%%%%%

%%%%%%%%%%%%%%%%%%%%%%%%%%%%%%%%%%%%%%%%%%%%%%%%%%%%%%%%%%%%%%%%%%

\end{document}